\RequirePackage{fix-cm}
\documentclass[smallextended]{svjour3}       % onecolumn (second format)
\smartqed  % flush right qed marks, e.g. at end of proof
\usepackage{graphicx}
\usepackage{textcomp}
\usepackage{epsfig}
\usepackage{amsmath,mathtools}
\usepackage{amssymb}
\usepackage[cp1250]{inputenc}
\usepackage{color}
\usepackage{esvect}
\usepackage{wasysym}

\usepackage{authblk}
\newcommand{\beq}{\begin{equation}}
\newcommand{\eeq}{\end{equation}}
\newcommand{\vep}{\varepsilon}
\newcommand{\beqa}{\begin{eqnarray}}
\newcommand{\eeqa}{\end{eqnarray}}
\newcommand{\bea}{\begin{eqnarray*}}
\newcommand{\ea}{\end{eqnarray*}}

\newcommand{\Lie}[1]{{\mathcal L}_{#1}}
\newcommand{\e}{{\mathrm{e}}}
\newcommand{\R}{\mathbb{R}}
\newcommand{\tran}{\mathsf{T}}
\newcommand{\X}{{X}}
\newcommand{\Y}{{Y}}
\usepackage[final,colorlinks]{hyperref}
\usepackage{url}
\definecolor{darkblue}{cmyk}{1,0,0,0.8}
\definecolor{darkred}{cmyk}{0,1,0,0.7}
\hypersetup{anchorcolor=black,
  citecolor=darkblue, filecolor=darkblue,
  menucolor=darkblue,pagecolor=darkblue,urlcolor=darkblue,linkcolor=darkblue}
\newtheorem{assumption}[lemma]{Assumption}
\DeclareMathOperator{\re}{Re}
\DeclareMathOperator{\spec}{spec}
\renewcommand{\d}{\mathrm{d}}
\begin{document}

\title{Dynamics of singularly perturbed sliding flow in Filippov systems%\thanks{Grants or other notes
%about the article that should go on the front page should be
%placed here. General acknowledgments should be placed at the end of the article.}
}
%\subtitle{Do you have a subtitle?\\ If so, write it here}

%\titlerunning{Short form of title}        % if too long for running head

\author{Piotr Kowalczyk         \and
        Jan Sieber %etc.
}

%\authorrunning{Short form of author list} % if too long for running head

\institute{Piotr Kowalczyk \at
              Department of Mathematics, Wroc{\l}aw University of Science and Technology, Poland \\
              \email{piotr.s.kowalczyk@pwr.edu.pl}           %  \\
%             \emph{Present address:} of F. Author  %  if needed
           \and
           Jan Sieber \at
              Department of Mathematics and Statistics, University of Exeter, United Kingdom \\
               \email{j.sieber@exeter.ac.uk}
}

\date{Received: date / Accepted: date}
% The correct dates will be entered by the editor

\maketitle

\begin{abstract}
In this article, we present an analysis of the effects of singular perturbations on the sliding motion in Filippov systems.
We show that singular perturbations may lead to qualitatively distinct topologies of phase space on the switching manifold, which we classify into six distinct topologies. Five of these topologies imply that singularly perturbed trajectory includes a segment (or segments) of sliding, and one topology, which we study here, implies the evolution characterised by switchings between trajectory segments along the switching surface, but without any sliding. In particular, we show that in the case of $n$-dimensional Filippov systems with one switching surface and $m$-dimensional fast dynamics, which plays the role of a stable singular perturbation, the flow follows sliding motion of the reduced system, but the perturbation becomes time dependent and may create a micro chaotic behaviour. However, a significant change in the flow direction is not possible. In the particular case of $1$-dimensional fast dynamics, the singular perturbation implies regular perturbation of $\mathcal{O}(\vep)$ of the perturbed flow at the points of switching between vector fields. 
% \PACS{PACS code1 \and PACS code2 \and more}
\end{abstract}

\section{Introduction}
\label{intro}

In the process of mathematical modelling, the standard procedure of reducing a given problem's complexity may result in a mathematical framework with discontinuous nonlinearities. To give just a few examples, systems characterised by impacts, switchings, and/or by events which on the macroscopic scale occur almost instantaneously \cite{AdDaNo:01,LeNi:04,diBeBuChaKo:08,MaLa:12,SiGl:24,Si:25}
%{coothwed2012} 
are all examples where the framework of discontinuous systems is of particular relevance. To such a class of systems belong so-called \emph{Filippov systems} \cite{Fi:88,Ut:92}. That is dynamical systems with phase space partitioned into non-empty adjacent regions where in each region a distinct set of smooth ordinary differential equations (ODEs) governs the dynamics. In the simplest case we me have just two regions, say $\mathcal{H}_+$ and $\mathcal{H}_-$ with a boundary $\mathcal{H}$ between these regions, called a \emph{switching manifold} (or \emph{surface}). We may then consider Filippov systems, differential equations of the form 
\begin{equation}
 \dot x = \begin{cases} f_+(x)\quad\quad \mbox{if}\quad\quad h(x) > 0, \\
                       f_-(x)\quad\quad \mbox{if}\quad\quad h(x) < 0,
         \end{cases}
         \label{eq:fil_system}
\end{equation}
with a boundary $\mathcal{H} = \{h(x) = 0\}$, where $f_\pm$ are vector fields obtained from ODEs describing the dynamics in $\mathcal{H}_\pm$ and $h(x)$ is a scalar valued multivariable smooth function. A particular feature of Filippov systems is the possibility that in some non-empty region $\mathcal{H}_\mathrm{ss}$ of the boundary $\mathcal{H}$  both vector fields point towards $\mathcal{H}$. In this case the directional derivatives satisfy
\begin{align*}
\frac{\d h(\X_+^t(x))}{\d t} = h_xf_+ < 0,\quad  \frac{\d h(\X_-^t(x))}{\d t} = h_x f_- > 0\mbox{\quad for $x\in\mathcal{H}_\mathrm{ss}\subseteq\mathcal{H}$,}
\end{align*}
where $\X_\pm^t$ are flows generated by the vector fields $f_\pm$ (all functions are taken at their argument $x$, and the subscript denotes the partial derivative). Hence, for points $x\in\mathcal{H}_\mathrm{ss}$ one can not define a flow of \eqref{eq:fil_system} by either following $\X_+^t$ or $\X_-^t$.
A convention intoduced by \cite{Fi:88,Ut:92} is to consider a \emph{sliding vector field}, which generates a \emph{sliding flow} $\X_\mathrm{s}^t$ on on $\mathcal{H}_\mathrm{ss}$. The sliding vector field is a convex combination of $f_\pm$ of the form
\begin{equation}
\begin{aligned}
\dot x &= f_s(x) = f_+ + \alpha(f_- - f_+)\mbox{,\ where} &
\alpha &= -\frac{h_xf_+}{h_x(f_- - f_+)}\mbox{, $\alpha \in [0,1]$.}
\end{aligned}
\label{eq:filippov_conv}
\end{equation}
%$$
%\alpha = -%\frac{h_xf_%+}{h_x(f_- %- %f_+)},\quad% \alpha \in [0, \, %1].
%$$
As we have already indicated in \cite{SiKo:10}, an important question to consider is how the sliding vector field is affected by perturbations. %As we have also indicated in \cite{SiKo:10} 
Regular (that is, $C^1$-small) perturbations applied to $f_\pm$ and $h$ do not qualitatively alter the system's dynamics provided some generic hyperbolicity and and transversality conditions hold, similar to the conditions and theorems on the persistence and stability of invariant sets in smooth dynamical systems \cite{fe:79}. An open question is how additional degrees of freedom that are stable and rapidly decaying  affect the sliding flow. Simple examples of such perturbations in modelling of systems with discontinuities are when one replaces a contact to a rigid wall by a finite stiffness, mass and damping in mechanical models, or takes into account a previously neglected small capacitance of an electrical element in a circuit with a switch. Mathematically these additional degrees of freedom are called stable \emph{singular perturbations}. The reverse approach is typically used to reduce the system dimension of an ODE model. One replaces rapidly converging subsystems of an ODE model by their equilibrium value, which changes slowly (\emph{quasistatically}), depending on the slower subsystems of the ODE. 

It is easy to demonstrate that such stable singular perturbations can change a sliding vector field into a rapid sequence of switching back-and-forth between vector fields  (see \cite{SiKo:10} for minimal examples). Thus, no direct equivalent of the persistence results for smooth vector fields by \cite{fe:79} can possibly hold in general. 
%in a way which implies a significantly different dynamics between that of the sliding flow in the reduced system, and the corresponding perturbed evolution in a slow-fast singularly perturbed Filippov system. Do we have a situation akin to smooth dynamical systems, where, in general, phenomena that persist under regular perturbations also persist under stable singular perturbations?
Fridman \cite{Fr:02a,Fr:02b} showed that, if the switching function depends only on the slow variables, periodic orbits persist under singular perturbations, but they acquire a small boundary layer. If one  allows that the switching decision function depends on fast variables, this may lead to a period-adding cascade or small-scale chaos that is not present in the reduced system \cite{SiKo:10}.

The aim of this paper is to close an existing gap in the theory of singularly perturbed Filippov systems and provide an answer to the question what is the effect of stable singular perturbation on sliding motion in Filippov systems. 
Our exploration how the Filippov convention interacts with additional stable degrees of freedom is complementary to recent work by Kristiansen and Hogan \cite{KaKr:19,KrHo:15a,KrHo:15b}. Their work starts from the assumption that the discontinuity in the right-hand side of \eqref{eq:fil_system} is the singular limit (``$\vep=0$'') of an underlying smooth singularly perturbed system (``$\vep>0$''). Then they study the dynamics for small $\vep>0$ using blow-up techniques \cite{KrSz:01} to observe limits of the dynamics and see if these limits are consistent with the Filippov convention. 

In the current work, we exclude the possible effects of the presence of so-called $U$-singularity (also termed as two-fold singularity) in Filippov systems of interests, which, among other works, have been studied in \cite{COLOMBO20131}.  
 %We should mention here that singularly perturbed Filippov system have been considered in a different framework using regularisation of the Filippov vector field and then a blow-up technique to determine the dynamics across the switching surface .

The outline of the rest of the paper is as follows. Section~\ref{sec:spfs} introduces the type of singular perturbations to Filippov systems we study, contrasting the \emph{full system} ($\vep>0$) and the \emph{reduced system} ($\vep=0$). Section~\ref{sec:reslom} performs a rescaling, zooming into the neighborhood of a point near the switching surface of the reduced system, a coordinate transformation and a truncation to leading order in $\vep$, after which $\vep$ disappears as parameter. The main results are then presented in Sec.~\ref{sec:results}, and in Sec. 5 and 6 where they are encapsulated in Lemma~\ref{thm:rep:switch},  Theorem~\ref{thm:convergence} and Theorem~\ref{thm:stable:fixedpoint}. Section ~\ref{sec:examples} illustrates our findings with a numerical example. Finally Section~\ref{sec:conclusions} concludes the paper. 

\section{Filippov systems and fast subsystems --- main results}
\label{sec:spfs}
\paragraph{Slow-fast vector fields with discontinuity and Filippov convention} We consider ODEs for $x\in\R^n$ with a discontinuous right-hand side, which interact with a \emph{fast subsystem} (an ODE for $y\in\R^m$ on a faster time scale):
\begin{align}
 \dot x &= \begin{cases} f_+(x,y,\vep)& \mbox{if $h(x,y,\vep) > 0$,} \\
                       f_-(x,y,\vep)& \mbox{if $h(x,y,\vep) < 0$,}
         \end{cases}
         \label{eq:s1}
         &%\\
         \vep\dot y &= g(x,y,\vep), 
 %\label{eq:s2}
\end{align}
where $\vep$ is small, right-hand sides $f_\pm:\mathbb{R}^{n + m + 1}\mapsto\mathbb{R}^{n}$, $g:\mathbb{R}^{n + m + 1}\mapsto\mathbb{R}^{m}$, and \emph{switching function} $h:\mathbb{R}^{n + m + 1}\mapsto\mathbb{R}$ are sufficiently smooth in some region in $\mathbb{R}^{n + m + 1}$ and the dot symbol denotes differentiation with respect to $t$. We permit that on the zero level set of $h$ (the \emph{switching manifold}, a codimension-$1$ hypersurface), %the \emph{switching function} $h(\cdot,\cdot,\vep)$,
$$
\mathcal{H}_\vep = \{(x,y)\in\mathbb{R}^{n + m}: h(x,y,\vep) = 0\},
$$
$f_+(x,y,\vep)$ and $f_-(x,y,\vep)$ may differ, such that the right-hand side in \eqref{eq:s1} is discontinuous in $\mathcal{H}_\vep$. The variable $y$ is called the \emph{fast variable}, while $x$ is the \emph{slow variable}. %We call $\mathcal{H}_0$ the \emph{switching manifold}
%(a codimension-$1$ surface in the full state space $\mathbb{R}^{n+m}$)

%To see how the Filippov sliding vector field \eqref{eq:filippov_conv} looks like for our slow-fast Filippov system, we express it in the form of equation \eqref{eq:fil_system}:
%\begin{equation}
%  \begin{pmatrix}
%      \dot x \\ \dot y 
%  \end{pmatrix} = \begin{cases} F_+(x,y,\vep) =  \begin{pmatrix} f_+(x,y,\vep) \\ \vep^{-1}  g(x,y,\vep)\end{pmatrix}&\mbox{if $h(x,y,\vep) > 0$,} \\[3ex]
%    F_-(x,y,\vep) =   \begin{pmatrix} f_-(x,y,\vep) \\ \vep^{-1}  g(x,y,\vep)\end{pmatrix}&\mbox{if $h(x,y,\vep) < 0$.}
%         \end{cases}        
% \label{eq:slow_fast_fil}
%\end{equation}
%Thus, the sliding vector field for the full fast-slow system \eqref{eq:slow_fast_fil} is
%\begin{equation}
%\begin{aligned}
%   \begin{pmatrix}\dot x \\ \dot y \end{pmatrix} &= F_s = F_+ + \alpha (F_- - F_+),&& \mbox{where}&
%   \alpha &= -\frac{h_xf_+ + (1/\vep) h_yg}{h_x(f_-- f_+)},
%\end{aligned}
%\label{eq:sl_full}
%\end{equation}
%and $f_\pm$, $g$, $h_x$ and $h_y$ are taken in in the point $(x,y,\vep)$ where the subscripts $x$, $y$ $\vep$ denote the partial derivatives.
%\paragraph{Transversal stability of fast dynamics} 
%The fast dynamics is obtained by rescaling time ($t=t_\mathrm{slow}=\vep t_\mathrm{fast}$), then setting $\vep$ to $0$ (such that $x$ is constant in time), and then determining the dynamics of $y$ depending on $x$. 
We assume that there exists a stable equilibrium $y_\mathrm{c}(x)$ ($y_\mathrm{c}:\R^n\to\R^m$) of the fast subsystem after rescaling time by $\vep$, $\dot y=g(x,y,\vep)$ at $\vep=0$ for each fixed $x$. %Let us define the map $y_\mathrm{c}:\mathbb{R}^{n} \to \mathbb{R}^{m}$ that assigns the equilibrium value for each $x$. 
The map $y_\mathrm{c}$  is defined implicitly by the equation
\begin{align}\label{ass:fast:equilibrium}
    g(x,y_\mathrm{c}(x), 0) = 0.
\end{align} 
The set of equilibria of the fast subsystem in rescaled time is called the \emph{slow manifold} 
\begin{align*}
\mathcal{M}_0 = \{(x,y)\in\mathbb{R}^{n+m}: g(x,y,0) = 0\}=\{(x,y)\in\mathbb{R}^{n+m}: y = y_\mathrm{c}(x) \},
\end{align*}
%is called 
%(a dimension-$n$, codimension $m$ manifold in $\mathbb{R}^{n+m}$).
\begin{assumption}[Transversal stability]\label{ass:fast:stable}
We assume that the equilibria in $\mathcal{M}_0$ are uniformly exponentially stable for fixed $x$ \textup{(}using $\spec A$ for the spectrum of a matrix $A$\textup{):}
% \begin{align}
%    \label{ass:fast:stable}
    $\re\spec g_y(x,y_\mathrm{c}(x),0) < -c_\mathrm{stab} < 0$ for all $x$.% such that $(x,y_\mathrm{c}(x),0)\in D$,}
%\end{align}    
\end{assumption}
We compare the dynamics of the full system \eqref{eq:s1} with that of the \emph{reduced system} on the slow manifold $\mathcal{M}_0$, which is also of Filippov type in $\R^n$,
\begin{align}
\label{ode:reduced}
    \dot x&=\begin{cases}
        f_+(x,y_\mathrm{c}(x),0)&\mbox{if $h_\mathrm{rd}(x)>0$,}\\
        f_-(x,y_\mathrm{c}(x),0)&\mbox{if $h_\mathrm{rd}(x)<0$,}
    \end{cases}\mbox{\quad where}&
    %\begin{aligned}
        h_\mathrm{rd}(x)&=\hspace*{1.2ex}h(x,y_\mathrm{c}(x),0).%\\
        %f_{\mathrm{rd},\pm}(x)&=f_\pm(x,y_\mathrm{c}(x),0).
    %\end{aligned}
\end{align}
\begin{assumption}[Attracting sliding in reduced system]\label{ass:red:sliding:stab}
    We assume that
\begin{align}
    \label{reduced:sliding:cond}
    h_{\mathrm{rd},x}f_+<0<h_{\mathrm{rd},x}f_-\mbox{,\quad where\quad}h_{\mathrm{rd},x}&=h_x-h_yg_y^{-1}g_x.
\end{align}
\end{assumption}
The coefficients in \eqref{reduced:sliding:cond} are evaluated at $(x,y_\mathrm{c}(x),0)$.
%\paragraph{Attracting sliding in reduced vector field}  with
The reduced system \eqref{ode:reduced} thus has an attracting sliding vector field in $\{x:h_{\mathrm{rd}}(x)=0\}$:
\begin{align}
    \label{reduced:sliding}
\begin{aligned}
    \dot x&=f_{\mathrm{rd,s}}(x):=f_++\alpha_\mathrm{rd}(f_--f_+)\mbox{, where\ }
    \alpha_\mathrm{rd}= -\frac{h_{\mathrm{rd},x}f_+}{h_{\mathrm{rd},x}(f_--f_+)}\mbox{,}
\end{aligned}
\end{align}
evaluating  all coefficients in $x$ and $(x,y_\mathrm{c}(x),0)$.
\begin{assumption}[Transversality of slow and switching manifolds]\label{ass:trans}\ \\
We assume that the intersection of the switching manifold  $\mathcal{H}_0$ and the slow manifold $\mathcal{M}_0$ is non-empty and transversal at $\vep=0$, such that
\begin{align*}
    \mathcal{H}_\mathrm{rd}=\mathcal{H}_0\cap\mathcal{M}_0&=\{(x,y)\in\R^{n+m}:h(x,y,0)=0, g(x,y,0)=0\}%\\
   % &=\{(x,y_\mathrm{c}(x))\in\R^{n+m}:h_{\mathrm{rd},x}(x)=h(x,y_\mathrm{c}(x),0)=0\}
\end{align*} 
is a smooth manifold of dimension $n-1$.
\end{assumption}
\paragraph{Main results}
We consider the dynamics of $(x,y)$ near a point $(x_0,y_\mathrm{c}(x_0))\in \mathcal{H}_\mathrm{rd}$ on the intersection of switching and slow manifolds. Which transitions between $f_-$, $f_+$ and sliding trajectories starting from $(x,y)$ are possible on time scale larger than $O(\vep)$ depends on the ordering of the four scalar quantities
\begin{align}\label{main:TS}
    T_\pm&:=h_x^0f_\pm^0,& S_\pm&:=h_y^0\left[g_y^0\right]^{-1}g_x^0f_\pm^0
\end{align}
(all coefficients evaluated at $(x_0,y_\mathrm{c}(x_0),0)$). When full system~\eqref{eq:s1} is restricted to half space $\{\pm h(x,y,\vep)>0\}$ it will have an invariant manifold that deviates from $y_\mathrm{c}(x)$ by order $\vep$. The quantities $S_\pm$ measure in which direction along $\mathcal{H}_0$ the intersection between $\mathcal{H}_0$ and $\mathcal{M}_0$ is split by this deviation. The quantities $T_\pm$ show where the full vector fields $(f_\pm,g/\vep)$ are tangent to the switching manifold $\mathcal{H}_0$, which can occur near $(x_0,y_\mathrm{c})$ despite Assumption~\ref{ass:red:sliding:stab} (note the presence of $h_{\mathrm{rd},x}$ instead of $h_x$ in \eqref{reduced:sliding:cond}). Assumption~\ref{ass:red:sliding:stab} is equivalent to $T_+<S_+$ and $S_-<T_-$. The six possible orderings of $S_\pm$ and $T_\pm$ are listed in Table~\ref{tab:tab1} and illustrated in Figure~\ref{fig_top1}. Of these two are symmetrically related, two permit only attracting sliding near $(x_0,y_\mathrm{c}(x_0))$, and one (the case $S_-<T_-<T_+<S_+$) permits only repelling sliding.
We then establish the following perturbation result for the case $S_-<T_-<T_+<S_+$ with repelling sliding.
\begin{theorem}[Smallness of stable fast perturbations]\label{thm:main:approx}
    Assume that the slow-fast system \eqref{eq:s1} satisfies Assumptions~\ref{ass:fast:stable}--\ref{ass:trans} near $(x_0,y_\mathrm{c}(x_0))$ and assume that $S_-<T_-<T_+<S_+$. There exists a perturbation bound $\rho_f>0$ such that for all $\rho_x,\rho_y,\rho_t>0$ we can find an $\vep_{\max}>0$ such that the following holds.
    
    If $|\delta_x|\leq\rho_x$, $|\delta_y|\leq\rho_y$, $t\in[0,\rho_t]$ and $\vep\in(0,\vep_{\max})$ then the trajectory  $(x(t),y(t))$ of \eqref{eq:s1} starting from $(x_0+\vep\delta_x,y_\mathrm{c}(x_0)+\vep\delta_y)$ without repelling sliding satisfies
    \begin{align}
        \label{main:eq:approx}
        \begin{bmatrix}
            x(\vep t)\\
            y(\vep t)    
        \end{bmatrix}-\begin{bmatrix}
            x(0)\\y(0)
        \end{bmatrix}&=\vep t
        \begin{bmatrix}
            f_\mathrm{rd,s}(x_0)\\
            \left[g_y^0\right]^{-1}g_x^0f_\mathrm{rd,s}(x_0)
        \end{bmatrix}+\vep r(t,\delta_x,\delta_y,\vep)        
    \end{align}
    for some continuous function $r$ bounded by 
$|r(t,x,y,\vep)|\leq \rho_f$ 
    for all $t\in\R$, $\vep\in[0,\vep_{\max}]$, $|x|\leq\rho_x$, $|y|\leq\rho_y$.
\end{theorem}
The theorem states that the reduced sliding flow \eqref{reduced:sliding} is a good approximation on time scales greater than $\vep$. This choice is possible in the theorem, since the bound $\rho_f$ of the perturbation scaled by $\vep$ is independent of the time scale $\rho_t$ scaled by $\vep$. Thus, we can choose $t\gg1$ such that still $\vep t\ll1$. In that case the drift of order $\vep t$ dominates the perturbation $\vep r$. Once $\vep t$ is no longer small, the drift term, which is evaluated at $(x_0,y_\mathrm{c}(x_0))$ is no longer approximately constant. The perturbation depends on $t$ and is only continuous so can be rapidly fluctuating. 
Informally, one could phrase the result of Theorem~\ref{thm:main:approx} as $x(t)$ for the full system \eqref{eq:s1} satisfying
\begin{align}
    \label{main:informal}
    \dot x(t)=f_\mathrm{rd,s}(x)+\tilde\vep \tilde r(t/\vep)+O(\vep),
\end{align}
where $\vep\ll\tilde\vep\ll1$ and $\tilde r$ is bounded but fluctuates rapidly (and will depend on quantities not shown in \eqref{main:informal}, such as $y(0)$). One obtains \eqref{main:informal} by dividing \eqref{main:eq:approx} by $\vep t$, assuming $t$ large and $\vep$ small, such that $\tilde\vep=1/t$.
In Section~\ref{sec:results} the perturbation is constructed as the trajectory of an approximately piecewise affine map for the scaled-up fast variable $\delta_y$. This map can be chaotic, such that we may have small ``micro chaos'' around the sliding flow, similar to what is observed when one digitizes control inputs \cite{haller1996micro,csernak2010digital,GlKo2010}. Small-scale fluctuations of the type in \eqref{main:eq:approx} were also constructed as a generic feature of \eqref{eq:s1} when it is near a grazing bifurcation \cite{SiKo:10}.

Our final main result, stated in Theorem \ref{thm:stable:fixedpoint} concerns the case when the equality in Theorem \ref{thm:main:approx} is exact, that is when the stable fast perturbation implies regular perturbation of $\mathcal{O}(\vep)$
for the slow-fast system (\ref{eq:s1}). This is the case when 
$$
            \begin{bmatrix}
            x(\vep t)\\    y(\vep t)
            \end{bmatrix}
            - 
            \begin{bmatrix}
            x(0)\\y(0)
            \end{bmatrix}
            =\vep t
        \begin{bmatrix}
            f_\mathrm{rd,s}(x_0)\\
            \left[g_y^0\right]^{-1}g_x^0f_\mathrm{rd,s}(x_0)
        \end{bmatrix}    
$$
at the switching points, with the correction term appearing as the regular perturbation of $\mathcal{O(\vep)}$, in the case $S_-<T_-<T_+<S_+$. As we show in Theorem \ref{thm:stable:fixedpoint} this is the situation when the fast dynamics $y$ is 1-dimensional. 
%\paragraph{Transversal intersection of switching and slow manifold}Furthermore,  Then the intersection set 
\section{Rescaling and reduction to leading order in small parameter}
\label{sec:reslom}
We will study the dynamics near a point $(x_0,y_\mathrm{c}(x_0))\in\mathcal{H}_\mathrm{rd}$, an arbitrary point on the switching manifold of the reduced system, where we have attracting sliding for the reduced system (so, \eqref{reduced:sliding:cond} is satisfied).
\subsection{Zoom-in in space and slow-down of time}
We now rescale the state variables 
\begin{align}
    x_\mathrm{old}& = x_0 + \vep x_\mathrm{new},& y_\mathrm{old} &= y_c(x_0) + \vep y_\mathrm{new},& t_\mathrm{old}&=\vep t_\mathrm{new}.
\label{eq:blow_up}
\end{align}
Substituting the zoom-in \eqref{eq:blow_up} into system equations \eqref{eq:s1}, we obtain (removing subscripts ``$\mathrm{new}$'' and ``$\mathrm{old}$'')
\begin{eqnarray}
\dot x & = &\begin{cases}f_+(x_0 + \vep x, y_c(x_0)+\vep y, \vep)& \mbox{if $h(x_0 + \vep x, y_c(x_0) + \vep y,\vep) > 0$,} 
\\
f_-(x_0 + \vep x, y_c(x_0)+\vep y, \vep)& \mbox{if $h(x_0 + \vep x, y_c(x_0) + \vep y,\vep) < 0$,}  
\end{cases}
\label{eq:sl}
\\ 
\vep\dot y & = & g(x_0 + \vep x, y_c(x_0) + \vep y,\vep). 
\label{eq:fst}
%\vep\dot u & = &\begin{cases}f_+(x_0 + \vep u, y_c(x_0)+\vep v, \vep)\quad \mbox{if}\quad h(x_0 + \vep u, y_c(x_0) + \vep v,\vep) > 0, \\
%f_-(x_0 + \vep u, y_0(x_0)+\vep v, \vep)\quad \mbox{if}\quad h(x_0 + \vep u, y_c(x_0) + \vep v,\vep) < 0,  
%\end{cases}
%\label{eq:sysrescls}
%\\ 
%\vep^2\dot v & = & g(x_0 + \vep u, y_c(x_0) + \vep v,\vep). 
%\label{eq:sysresclf}
\end{eqnarray} 
Note that using the above rescaling means that a distance of $\mathcal{O}(\vep)$ in $(x_\mathrm{old},y_\mathrm{old})$ from $(x_0,y_\mathrm{c}(x_0))$ corresponds to a distance of $\mathcal{O}(1)$ in $(x_\mathrm{new},y_\mathrm{new})$ from $(0,0)$. In the new coordinates $(x_\mathrm{new},y_\mathrm{new})=(0,0)$ corresponds to $(x_\mathrm{old},y_\mathrm{old})=(x_0,y_\mathrm{c}(x_0))$.

%We further rescale the time variable such that $t = \mathcal{O}(\vep)$ will correspond to $\tau  = \mathcal{O}(1)$. That is we introduce $\tau = (1/\vep) t$. 
%Substituting $\tau$ into system equations (\ref{eq:sysrescls})--(\ref{eq:sysresclf}) gives the system
%\begin{eqnarray}
%\dot u & = &\begin{cases}f_+(x_0 + \vep y, y_c(x_0)+\vep v, \vep)\quad \mbox{if}\quad h(x_0 + \vep u, y_c(x_0) + \vep v,\vep) > 0, 
%\\
%f_-(x_0 + \vep y, y_c(x_0)+\vep v, \vep)\quad \mbox{if}\quad h(x_0 + \vep u, y_c(x_0) + \vep v,\vep) < 0,  
%\end{cases}
%\label{eq:sl}
%\\ 
%\vep\dot v & = & g(x_0 + \vep u, y_c(x_0) + \vep v,\vep). 
%\label{eq:fst}
%\end{eqnarray} 
\subsection{Truncation to leading order in $\vep$ and shift of origin}
Expanding the full system \eqref{eq:sl}--\eqref{eq:fst} to leading order in $\vep$ we get
\begin{eqnarray}
\dot x & = &\begin{cases}f_+^0 + \mathcal{O}(\vep)& \mbox{if $h_x^0x + h_y^0y + h_\vep^0 + \mathcal{O}(\vep) > 0$,} 
\\
f_-^0+\mathcal{O}(\vep) & \mbox{if $h_x^0x + h_y^0y + h_\vep^0 + \mathcal{O}(\vep) < 0$,}  
\end{cases}
\label{eq:slexp}
\\ 
%\mbox{\jsq{$\vep$?}}
\dot y & = & g_x^0 x + g_y^0 y + g_\vep^0 + \mathcal{O}(\vep)\mbox{, \quad where}\\
\label{eq:fstexp}
\nonumber
f_\pm^0&=&f_\pm(x_0,y_\mathrm{c}(x_0),0),
\end{eqnarray} 
and $h_x^0\in \mathbb{R}^{1\times n}$,  $h_y^0\in \mathbb{R}^{1\times m}$, $h_\vep^0\in \mathbb{R}^{1\times 1}$, $g_x^0\in \mathbb{R}^{m\times n}$,  $g_y^0\in \mathbb{R}^{m\times m}$ and $g_\vep^0\in \mathbb{R}^{m\times 1}$ are the Jacobian matrices evaluated at $(x_0,y_c(x_0),0)$. Correspondingly, the quantity
\begin{align}
    h_{\mathrm{rd},x}^0&=h_x^0-h_y^0[g_y^0]^{-1}g_x^0,
    \label{eq:hred}
\end{align}
is the normal to the switching surface in $(x_0,y_\mathrm{c}(x_0))$ of the reduced system \eqref{ode:reduced}, which is non-zero as we assume that the expansion point is on the switching manifold of the reduced system by assumption \eqref{reduced:sliding:cond}.
%\jsq{where is this needed?} without loss of generality (w.l.o.g.), we assume that $h_x^0\not = 0$.
By shifting the rescaled $x$ and $y$ coordinates we may also eliminate the constant terms in the right-hand side expansion ($g_\vep^0$) and the switching function ($h_\vep^0$):\begin{align*}
    x_\mathrm{old}&=x_\mathrm{new}+\delta_x,&
    y_\mathrm{old}&=y_\mathrm{new}+\delta_y.   
\end{align*} 
The shifts $\delta_x$ and $\delta_y$ have to satisfy the affine system of $m+1$ equations
\begin{align*}
    0&=h_x^0\delta_x+h_y^0\delta_y+h_\vep^0,&
    0&=g_x^0\delta_x+g_y^0\delta_y+g_\vep^0,  
\end{align*}
for which one possible solution is
\begin{align*}
    \delta_x&=\delta_h h_{\mathrm{rd},x}^\tran,&\delta_y&=-\left[g_y^0\right]^{-1}\left[g_x^0\delta_hh_{\mathrm{rd},x}^\tran+g_\vep^0\right],\mbox{\ with\ } \delta_h=\frac{h_y^0[g_y^0]^{-1}g_\vep^0-h_\vep^0}{h_{\mathrm{rd},x}h_{\mathrm{rd},x}^\tran}
\end{align*}
(this choice of $\delta_x$ shifts orthogonal to the switching surface). In the new (zoomed-in, slowed-down, and shifted) coordinates, the full system of equation has the simplified form
%\paragraph{Shift in $v$}

%For the brevity of notation, let $$g_u^0 = B\mbox{,}\, g_v^0 = -A \mbox{,}\, g_\vep^0\mbox{,} = C\mbox{,}\,$$ and with some abuse of notation, in what follows, we drop the ``$0$'' superscript.
%We may shift $v$ by $\delta_v$  according to:
%$$
%A\delta_v = B\delta_u + C \Leftrightarrow \delta_v = A^{-1}B\delta_u + A^{-1}C,
%$$
%where we can isolate $\delta_v$ because $A$ is invertible by our assumption that the singular perturbation is stable. 
%This shift eliminates $g_\vep^0$ in equation (\ref{eq:fstexp}).

%We compute $\delta_u$ so that term $h_\vep^0$ is eliminated in equation (\ref{eq:slexp}). Namely, we have
%$$
%h_u \delta_u + h_vA^{-1}B\delta_u + h_vA^{-1}C = -h_\vep \Leftrightarrow (h_u+h_vA^{-1}B)\delta_u = -(h_\vep + h_vA^{-1}C).
%$$ 
%As we assume that point $(x_0,y_c(x_0))$ is on the attracting sliding part of the switching manifold in the reduced system we have that
%Let us assume that 
%$h_u+h_vA^{-1}B\not = 0$. Thus, we may choose a translation $\tilde u = u + \delta_u$ and $\tilde v = v + \delta_v$ such that our model system (\ref{eq:slexp})--(\ref{eq:fstexp}), to leading order, simplifies to   
\begin{eqnarray}
\dot x & = &\begin{cases}f_+^0+O(\vep) & \mbox{if $h_x^0 x + h_y^0 y +O(\vep)> 0$,} 
\\
f_-^0+O(\vep)&\mbox{if $h_x^0 x + h_y^0 y +O(\vep)< 0$,}  
\end{cases}
\label{eq:sllo:eps}
\\ 
\dot y & = & g_x^0x+g_y^0 y+O(\vep), 
\label{eq:fstlo:eps}
\end{eqnarray} 
where the two $O(\vep)$ terms in the switching condition are identical.

\paragraph{Truncated reduced system} The truncated reduced system is is the piecewise constant flow
\begin{align}
\dot x & = \begin{cases}f_+^0 & \mbox{if $h_{\mathrm{rd},x}^0x > 0$,} 
\\
f_-^0 & \mbox{if $h_{\mathrm{rd},x}^0x < 0$,} 
\end{cases}
\label{eq:rS}
\\
0 & =  g_x^0x+g_y^0 y\mbox{\quad (thus, $y=-[g_y^0]^{-1}g_x^0 x$).}
\label{eq:rSm}
\end{align}
In particular, we observe that the condition $h_{\mathrm{rd},x}^0f_+^0 < 0 < h_{\mathrm{rd},x}^0f_-^0 $ for attracting sliding for the truncated reduced flow  \eqref{eq:rS},\,\eqref{eq:rSm}, spelled out 
\begin{equation}
h_{\mathrm{rd},x}^0f_+^0=\left[h_x^0 - h_y^0\left[g_y^0\right]^{-1}\!\!g_x^0\right]f_+^0<0<\left[h_x^0 - h_y^0\left[g_y^0\right]^{-1}\!\!g_x^0\right]f_-^0=h_{\mathrm{rd},x}^0f_-^0,
\label{eq:slCrS}
\end{equation}
also implies sliding of the reduced flow of \eqref{eq:sllo:eps},\,\eqref{eq:fstlo:eps} for all small $\vep$ and all $x,y$ of order $1$. 
%Since in what follows all the constant quantities are evaluated at $(0,0)$ we omit the superscript ``$0$'' in $h_{\mathrm{rd},x}^0$, $f_+^0$ and all the other related quantities. 
%We dropped the tilde symbol for the clarity of notation in equation (\ref{eq:fstlo}).
%\subsection{Reduced system}
%On the set $\{v = A^{-1}Bu \}$, $\dot v = 0$, and in such case the system dynamics is determined by the dynamics of the vector field
%Sliding condition is then given by:
%\subsection{Slow-fast system}
%We are interested how the sliding flow of the reduced system (\ref{eq:rS}) is affected by fast dynamics in the full system. 
%\subsection{The sliding vector field of the truncated reduced system}
As we assume that the sliding condition \eqref{eq:slCrS} holds, \eqref{eq:rS} has a constant sliding vector field
on the set
\begin{align*}
    \mathcal{H}_{\mathrm{rd},x}:=\{h_{\mathrm{rd},x}^0 x = 0 \}, 
\end{align*} 
(the $x$-components of the discontinuity set $\mathcal{H}_\mathrm{rd}$ in the new coordinates) of the form
\begin{align}
\dot x  = f_\mathrm{rd,s}^0=(1 - \alpha_\mathrm{rd}^0)f_+^0 + \alpha_\mathrm{rd}^0 f_-^0=
f_+^0 + \alpha_\mathrm{rd}^0(f_-^0 - f_+^0),
\label{eq:slRf}
\end{align}
where $\alpha_\mathrm{rd}^0 \in [0,1] $ is determined by the condition that $f_\mathrm{rd,s}^0$ is in $\mathcal{H}_{\mathrm{rd},x}$, satisfying $h_{\mathrm{rd},x}^0 f_\mathrm{rd,s}^0 = 0$:
\begin{equation}
\alpha_\mathrm{rd}^0 = \frac{h_{\mathrm{rd},x}^0 f_+^0}{h_{\mathrm{rd},x}^0(f_+^0 - f_-^0)}.
\label{eq:alpha}
\end{equation}
So the constant sliding vector field $f_\mathrm{rd,s}^0$ for the truncated reduced system \eqref{eq:rS} equals
\begin{align}
f_\mathrm{rd,s}^0 &= f_+^0 -\frac{h_{\mathrm{rd},x}^0 f_+^0}{h_{\mathrm{rd},x}^0(f_-^0 - f_+^0)}(f_-^0 - f_+^0)=\frac{f_+^0h_{\mathrm{rd},x}^0f_-^0 - f_-^0h_{\mathrm{rd},x}^0f_+^0}{h_{\mathrm{rd},x}^0(f_-^0-f_+^0)}
\label{eq:slRf3}
\\
&=\frac{f_+^0(h_x^0f_-^0 - h_y^0\left[g_y^0\right]^{-1}g_x^0f_-^0) - f_-^0(h_x^0f_+^0 - h_y^0\left[g_y^0\right]^{-1}g_x^0f_+^0)}{h_x^0(f_-^0-f_+^0) - h_y^0\left[g_y^0\right]^{-1}g_x^0(f_-^0 -f_+^0)}.
\label{eq:slRf4}
\end{align}
The reduced system after truncation, \eqref{eq:slRf}, has, thus, a parallel flow in direction $f_\mathrm{rd,s}^0$.

\paragraph{Truncated full system} Neglecting terms of order $\vep$ the truncated full system is
\begin{align}
%\label{eq:sllo}
\dot x & = \begin{cases}f_+^0 & \mbox{if $h_x^0 x + h_y^0 y> 0$,} 
\\ 
f_-^0&\mbox{if $h_x^0 x + h_y^0 y< 0$,}  
\end{cases}
& 
\label{eq:fstlo}
\dot y & =  g_x^0x+g_y^0 y. 
\end{align} 
We change the $y$ coordinate such that it moves with $x$, eliminating the $x$ dependence in the $\dot y$ part of \eqref{eq:fstlo}: $y_\mathrm{new}= -g_y^0y_\mathrm{old}-g_x^0x$. In new coordinates (dropping subscript ``new'') the full truncated system satisfies
\begin{equation}
\begin{bmatrix}
\dot x \\
\dot y 
\end{bmatrix} = 
\begin{cases}
F_+ = 
\begin{bmatrix}
f_+^0 \\
 g_y^0y-g_x^0f_+^0
\end{bmatrix}&\mbox{if $h_{\mathrm{rd},x}^0x - h_y^0\left[g_y^0\right]^{-1}y > 0$,} \\[2.5ex]
F_- = \begin{bmatrix} f_-^0 \\
 g_y^0y-g_x^0f_-^0
\end{bmatrix}& \mbox{if $h_{\mathrm{rd},x}^0x - h_y^0\left[g_y^0\right]^{-1}y < 0$,}
\end{cases}
\label{eq:slfxy}
\end{equation}
and the discontinuity manifold in these new coordinates is the subspace
\begin{align*}
\mathcal{H}_0=\left\{(x,y):h_{\mathrm{rd},x}^0x = h_y^0\left[g_y^0\right]^{-1}y\right\}.
\end{align*}
As long as the full truncated system \eqref{eq:slfxy} stays in the respective domain for $F_\pm$, the affine subspace
\begin{align}\label{eq:slowman}
        \mathcal{M}_{\mathrm{slow},\pm}&=\{(x,y): y=y_{\mathrm{sl},\pm}\mbox{, $x$ arbitrary}\}\mbox{,  where}&y_{\mathrm{sl},\pm}&:=\left[g_y^0\right]^{-1}g_x^0f_\pm^0,
\end{align}
is exponentially attracting. Condition \eqref{eq:slCrS}, that the reduced system has attracting sliding ($h_{\mathrm{rd},x}f_+^0<0<h_{\mathrm{rd},x}f_-^0$) is equivalent to the conditions
\begin{align}
    h_x^0f_+^0&<h_y^0y_{\mathrm{sl},+},&
h_x^0f_-^0&>h_y^0y_{\mathrm{sl},-}
\end{align}
on $y_{\mathrm{sl},\pm}$.
%We will comment for each result for \eqref{eq:slfxy} on how it persists for small $\vep$. 
Assuming that \begin{align}
    \label{eq:full:gen} h_x^0(f_+^0-f_-^0)\neq 0,
\end{align}the full truncated system has 
\begin{itemize}
    \item attractive sliding if $h_x^0f_+^0<h_x^0f_-^0$ for $(x,y)\in\mathcal{H}_0$ where $y$ satisfies
\begin{align}\label{full:attr:sliding}
    h_x^0f_+^0&<h_y^0y<h_x^0f_-^0
\end{align}
\item repelling sliding  if $h_x^0f_+^0>h_x^0f_-^0$ for $(x,y)\in\mathcal{H}_0$ where $y$ satisfies
\begin{align}\label{full:rep:sliding}
    h_x^0f_+^0&>h_y^0y>h_x^0f_-^0
\end{align}
\end{itemize}
(conditions are independent of the $x$ component). The boundary of these sliding regions are the affine tangent subpaces
 \begin{equation}
 \mathcal{H}_{\mathrm{tan},\pm}=\{(x,y)\in\mathcal{H}_0:h_y^0y=h_x^0f_\pm^0\}
% , 
%h_y^0 y_+^{(T)} = h_x^0f_+^0\quad \mbox{and}\quad h_y^0y_-^{(T)} = h_x^0f_-^0\mbox{.}
\label{eq:tanpoints}
 \end{equation}
(these are subspaces of the switching subspace $\mathcal{H}_0$). The locations of the subspaces $\mathcal{M}_{\mathrm{slow},\pm}$ and $\mathcal{H}_{\mathrm{tan},\pm}$ relative to each other will determine the possible signatures of the dynamics (sequencing switches and sliding along trajectories) in the full truncated system in Section~\ref{sec:results}.

For sliding in the full truncated system \eqref{eq:slfxy} the factor $\alpha$ for the Filippov convention \eqref{eq:filippov_conv} simplifies to (assuming $h_x^0f_+^0\neq h_x^0f_-^0$ as in \eqref{eq:full:gen})
\begin{equation}
\alpha = \frac{h_x^0f_+^0 - h_y^0 y}{h_x^0(f_+^0-f_-^0)}, 
\label{eq:alpa_fs}
\end{equation}
which depends on $y$. Thus, the sliding vector field $F_\mathrm{s}$ of the truncated full system \eqref{eq:slfxy} is an affine flow:
\begin{align}
\label{eq:slidingFull}
\begin{bmatrix}
\dot x \\ \dot y 
\end{bmatrix} &= %&
F_\mathrm{s}  =b_\mathrm{s}+A_\mathrm{s}
\begin{bmatrix}
x \\ y
\end{bmatrix}\mbox{,\ where}\\
\nonumber
b_\mathrm{s} &= \begin{bmatrix}
f_+^0 \\
 -g_x^0f_+^0
\end{bmatrix} 
-\frac{h_x^0f_+^0}{h_x^0(f_-^0-f_+^0)} \begin{bmatrix}
f_-^0 - f_+^0\\
 -g_x^0(f_-^0 - f_+^0), &
\end{bmatrix}
\mbox{, and\ }\\ 
\nonumber
A_\mathrm{s}& = 
\begin{bmatrix}
0 & & \cfrac{(f_-^0 - f_+^0)h_y^0}{h_x^0(f_-^0 - f_+^0)} \\
 0& & g_y^0 - \cfrac{g_x^0(f_-^0 - f_+^0)h_y^0}{h_x^0(f_-^0 - f_+^0)}
\end{bmatrix}\mbox{.}
%\\
%\begin{bmatrix}
%f_+^0 \\
% g_y^0y-g_x^0f_+^0
%\end{bmatrix}
%-\frac{h_x^0f_+^0 - h_y^0 y}{h_x^0(f_-^0-f_+^0)} \begin{bmatrix}
%f_-^0 - f_+^0\\
% -g_x^0(f_-^0 - f_+^0)
%\end{bmatrix} \\ \nonumber
%& 
%=  
%\begin{bmatrix}
%f_+^0 \\
% -g_x^0f_+^0
%\end{bmatrix} 
%-\frac{h_x^0f_+^0}{h_x^0(f_-^0-f_+^0)} \begin{bmatrix}
%f_-^0 - f_+^0\\
% -g_x^0(f_-^0 - f_+^0)
%\end{bmatrix} %\\
%&  & 
%+
%\begin{bmatrix*}[r]
%0 & \cfrac{(f_-^0 - f_+^0)h_y^0}{h_x^0(f_-^0 - f_+^0)} \\[2ex]
% 0 & g_y^0 - \cfrac{g_x^0(f_-^0 - f_+^0)h_y^0}{ h_x^0(f_-^0 - f_+^0)}
%\end{bmatrix*}
%\begin{bmatrix}
%x \\ y
%\end{bmatrix}\\ %\nonumber 
%& 
%= C + A \left(\begin{array}{c}
%x \\ y
%\end{array}\right)\mbox{,}
\end{align} 
%where 
%$$
%C = \begin{bmatrix}
%f_+^0 \\
% -g_x^0f_+^0
%\end{bmatrix} 
%-\frac{h_x^0f_+^0}{h_x^0(f_-^0-f_+^0)} \begin{bmatrix}
%f_-^0 - f_+^0\\
% -g_x^0(f_-^0 - f_+^0)
%\end{bmatrix}\mbox{, and\ } A = 
%\begin{bmatrix}
%0 & \cfrac{(f_-^0 - f_+^0)h_y^0}{h_x^0(f_-^0 - f_+^0)} \\
% 0 & g_y^0 - \cfrac{g_x^0(f_-^0 - f_+^0)h_y^0}{h_x^0(f_-^0- f_+^0)}
%\end{bmatrix}\mbox{.}
%$$
We observe that the matrix $A_\mathrm{s}$ in the linear term of $F_\mathrm{s}$ has the $n$-dimensional nullspace $\{(x,y):y=0, x \mbox{\ arbitrary}\}$. Furthermore, the dynamics of $F_\mathrm{s}$ does not need to be linearly stable even though $g_y^0$ is assumed to generate stable dynamics because of the non-small rank-$1$ perturbation $g_x^0(f_-^0 - f_+^0)h_y^0/(h_x^0(f_-^0 - f_+^0))$.
The two conditions, \eqref{eq:slCrS} and \eqref{eq:full:gen},
\begin{align*}
    h_x^0(f_-^0-f_+^0)&\neq 0,&
    h_{\mathrm{rd},x}^0(f_-^0 - f+^0&)> 0
\end{align*}
%Together with the condition on the existence of attracting sliding in the reduced system \eqref{eq:slCrS}, $ h_{\mathrm{rd},x}^0(f_-^0 - f+^0) = h_x^0(f_-^0 - f+^0) - h_y^0\left[g_y^0\right]^{-1}g_x^0(f_-^0-f_+^0) > 0$ these two conditions 
ensure that the U-singularity \cite{diBeBuChaKo:08} occurs neither in the sliding vector field $f_\mathrm{rd,s}^0$ of the reduced truncated system, nor in the affine sliding vector field $F_\mathrm{s}$ of the full truncated system.  

\section{Classification of switching signatures}
\label{sec:results}

\begin{figure}
\centering{
\includegraphics[width=1\textwidth]{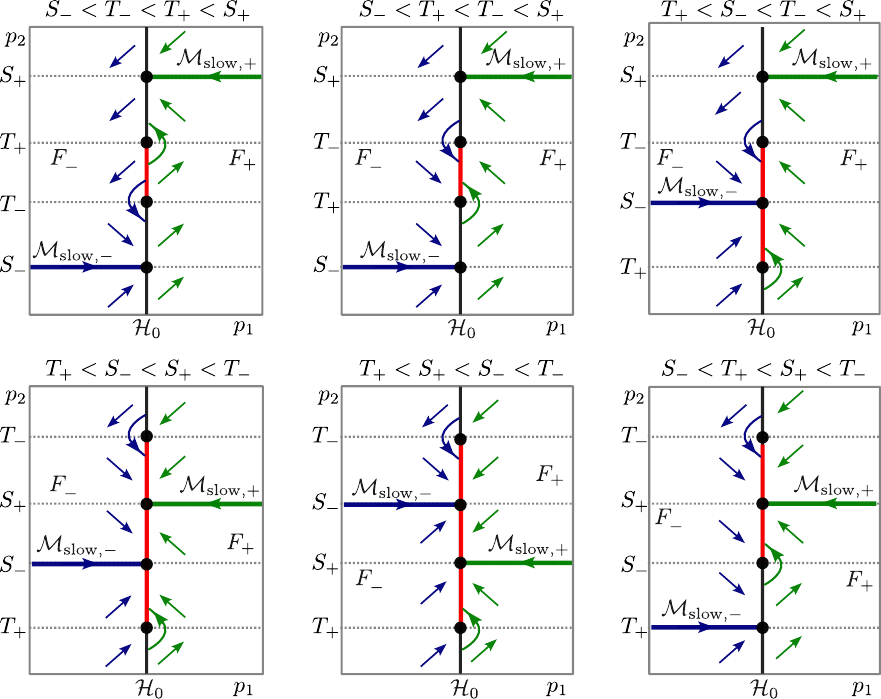}
}
\caption{Projection of phase portraits into plane by projection $(p_1,p_2)=P(x,y)=(h_{\mathrm{rd},x}^0x-h_y^0\left[g_y^0\right]^{-1}y,h_y^0y)$ for cases listed in Table~\ref{tab:tab1}.  Note that arrows may cross if $\dim y=m>1$. The sign of the horizontal component of all arrows is correctly reflected in the sketch. The vertical component is only illustrative (indicating attractivity of $\mathcal{M}_{\mathrm{slow},\pm}$ for $F_\pm$). Sections~\ref{sec:swbeh} and \ref{sec:1dfast} study the case in the top-left panel, $S_-<T_-<T_+<S_+$.%depiction of the topology of the switching surface and the relation to vector fields on either side of the switching surface corresponding to Cases 1 to 6 in Tab. ~\ref{tab:tab1}. 
}
\label{fig_top1}
\end{figure}

 The different possible switching signatures can be classified by considering the relative locations of the attracting subspace $\mathcal{M}_{\mathrm{slow},\pm}$ for each of the vector fields  $F_\pm$ near the switching subspace $\mathcal{H}_0$ and the tangency subspaces $\mathcal{H}_{\mathrm{tan},\pm}$, which determine where sliding regions in $\mathcal{H}_0$ start or end. Conveniently, all subspaces are determined by the projection $h_y^0$ of the $y$-component of the phase space variable, such that for each $(x,y)$ we may consider its linear projection into the plane
 \begin{align*}
    P&:\R^{n+m}\to \R^2\mbox{\quad with}\\
    [P_1&(x,y),P_2(x,y)]=\left[h_{\mathrm{rd},x}^0x-h_y^0\left[g_y^0\right]^{-1}y,\ h_y^0y\right].    
 \end{align*}
 In projection $P$ the switching subspace $\mathcal{H}_0$ is the vertical axis
 \begin{align*}
 \mathcal{H}_0=\{(x,y): P_1(x,y)=0\},    
 \end{align*} such that the vector fields $F_\pm$ are applicable in the half spaces projecting onto the positive and negative half-planes \begin{align*}
    \mathcal{H}_\pm=\{(x,y):\pm\, P_1(x,y)>0\}.
\end{align*}
The  projections of the tangent subspaces for sliding in $\mathcal{H}_0$ are the points 
 \begin{align*}
     &(0,T_\pm)\mbox{,\quad where}&T_\pm&=h_x^0f_\pm^0\mbox{,\quad such that}&\mathcal{H}_{\mathrm{tan},\pm}&=\{(x,y):P_2(x,y)=T_\pm\}.
\end{align*}
The projections of the attracting invariant subspaces $\mathcal{M}_{\mathrm{slow},\pm}$ of $F_\pm$ for the regions $\mathcal{H}_\pm$ are the horizontal levels
\begin{align*}
     \mathcal{M}_{P,\mathrm{slow},\pm}&=\left\{(p_1,p_2): \pm p_1>0, p_2=S_\pm\right\}\mbox{\ where}\\ S_\pm&=h_y^0y_{\mathrm{sl},\pm}\mbox{, such that}\\
     \mathcal{M}_{\mathrm{slow},\pm}&=\{(x,y):\pm P_1(x,y)>0, P_2(x,y)=S_\pm\}.
\end{align*}
In projection $P=(p_1,p_2)$, the nullclines for the horizontal component of motion in half plane $\{\pm\,p_1>0\}$ are the horizontal levels $T_\pm$: $\{\pm\,p_1>0, p_2=h_y^0f_\pm^0\}$. The vertical component of motion is attracted toward $\mathcal{M}_{P,\mathrm{slow},\pm}$ in each half-plane. If $y$ is scalar (such that we have only one fast variable in the original system \eqref{eq:s1}), the line $\mathcal{M}_{P,\mathrm{slow},\pm}$ acts as a nullcline for vertical motion. For $\dim y=m>1$ this may not be the case such that the vertical component of arrows in Figure~\ref{fig_top1} is for illustration only.
By our assumption \eqref{eq:slCrS} that the reduced flow has attractive sliding the levels $T_\pm$ and levels $S_\pm$ must satisfy
\begin{align}
    \label{eq:tpm:spm}
    T_-&>S_-,&
    S_+&>T_+,
\end{align}
resulting in the six possible arrangements of $S_\pm$ and $T_\pm$ illustrated in Figure~\ref{fig_top1} and classified in Table~\ref{tab:tab1}.
\begin{table}[t]
    \centering
    \begin{tabular}{l@{\hspace*{2em}}l@{\hspace*{2em}}l}
    \hline\\[-1.5ex]
    Vertical ordering& Sliding & Possible long-term switching\\[1ex]
    $S_-<T_-<T_+<S_+$ & repelling & $+-$, $-+$\\
    $S_-<T_+<T_-<S_+$ & attracting & $+-$, $-+$, $\mathrm{s}$, $\mathrm{s}+$, $\mathrm{s}-$, $+\mathrm{s}$, $-\mathrm{s}$\\
    $T_+<S_-<T_-<S_+$ & attracting & $-\mathrm{s}$,$\mathrm{s}-$, $\mathrm{s}$\\
    $T_+<S_-<S_+<T_-$ & attracting & $\mathrm{s}$ \\ 
    $T_+<S_+<S_-<T_-$ & attracting & $\mathrm{s}$\\
    $S_-<T_+<S_+<T_-$ & attracting & $+\mathrm{s}$,$\mathrm{s}+$, $\mathrm{s}$\\[1ex]
    \hline % \hline   
    \end{tabular}
    \caption{Possible transitions of dynamics near switching subspace $\mathcal{H}_0$ of truncated full systems \eqref{eq:slfxy}. Left column: ordering of $T_\pm$ for tangent set $\mathcal{H}_{\mathrm{tan},\pm}$ and of $S_\pm$ for attracting subspace $\mathcal{M}_{P,\mathrm{slow},\pm}$ (see vertical ordering of $S_\pm$ and $T_\pm$ in Figure~\ref{fig_top1}). Middle column: nature of sliding on sliding segment (attracting or repelling) according to \eqref{full:attr:sliding} and \eqref{full:rep:sliding}. Right column: switches (between $\mathcal{H}_+$, $\mathcal{H}_-$ and sliding that are possible in the long term: symbol $+-$ stands for crossing from $\mathcal{H}_+$ to $\mathcal{H}_-$, $-+$ vice versa, $\mathrm{s}\pm$ for leaving sliding toward $\mathcal{H}_\pm$, $\pm\mathrm{s}$ for entering sliding from $\mathcal{H}_\pm$, single $\mathrm{s}$ for long-term sliding. See also Figure~\ref{fig_top1} for illustration.}
    \label{tab:tab1}
\end{table}
The cases $T_+<S_-<T_-<S_+$ and $S_-<T_+<S_+<T_-$ are equivalent when interchanging the symbols $+$ and $-$ and the vector fields $F_+$ and $F_-$. Thus, overall, we have five qualitatively different scenarios up to change of subscripts $+$ and $-$. 
As Figure~\ref{fig_top1} shows, the sliding region along the vertical axis $\{p_1=0\}$ always occurs between points $T_+$ and $T_-$ (and sliding is attracting if $T_+<T_-$). The right column in Table~\ref{tab:tab1} shows switching transitions that are possible in the long-term dynamics near $\mathcal{H}_0\cap(\mathcal{M}_{P,\mathrm{slow},+}\cup \mathcal{M}_{P,\mathrm{slow},-})$. We observe that there is an upper limit for the time the flow can follow $F_\pm$ such that every period of following $F_\pm$ (indicated by symbol $\pm$ in Table~\ref{tab:tab1}) must be followed by a switch, either to $F_\mp$ (symbol $\mp$ in Table~\ref{tab:tab1}), or to $F_\mathrm{s}$ (symbol $\mathrm{s}$ in Table~\ref{tab:tab1}).
     
\paragraph{Scenario $S_-<T_-<T_+<S_+$ (top row in Table~\ref{tab:tab1}, see Figure~\ref{fig_top1})}
As the sliding set is repelling %The first scenario we describe is depicted in Fig.~\ref{fig_top1} Case 1. The arrangement of sets $S_\pm$, $T_\pm$ is given by $S_+ > T_+ > T_- > S_-$.
between the levels $T_+$ and $T_-$ there exists a repelling sliding set. This implies that any trajectory will follow $F_+$ and $F_-$ alternatingly with infinitely many switches from $F_-$ to $F_+$ between $S_-$ and $T_-$, and from $F_+$ to $F_-$ between $T_+$ and $S_+$. In Figure~\ref{fig_top1} this is schematically indicated by the arrows pointing away from the vertical line depicting the switching manifold. These arrows indicate the directions of the vector fields $F_\pm$. %Between sets on the switching surface corresponding to $S_+$ and $T_+$, and  $S_-$ and $T_-$, there exist sets where the vector fields $F_\pm$ point in the same direction with respect to the switching surface. This implies that starting an evolution on the switching surface, say by following $F_+$, the trajectory will return to it and switch to $F_-$ and then again to $F_+$, and so on. We will explain this in detail in Sec. \ref{sec:swbeh}. 
\paragraph{Scenario $S_-<T_+<T_-<S_+$ (row $2$ in Table~\ref{tab:tab1}, see Figure~\ref{fig_top1})}
%The second scenario we describe is depicted in Fig.~\ref{fig_top1} Case 2. Note that the sets $T_+$ and $T_-$ swap sides when compared with Scenario 1 and so the arrangement of sets $S_\pm$, $T_\pm$ is given by $S_+ > T_- > T_+ > S_-$. 
In the region between $T_+$ and $T_-$ there exists an attracting sliding set. This is schematically indicated by the arrows pointing towards the vertical line depicting the switching manifold.  Every transition between flows is possible in the long run (apart from permanently following $F_+$ or permanently following $F_-$ not being possible). Which sequence actually occur depends on the coefficients of the problem. %Similarly as before, between sets on the switching surface corresponding to $S_+$ and $T_-$, and  $S_-$ and $T_+$, there exist sets where the vector fields $F_\pm$ point in the same direction with respect to the switching surface. This implies that starting an evolution on the switching surface, we may follow the sliding vector field. If the sliding vector field is stable we are confined to it. Otherwise sets $T_-$ or $T_+$ may be reached, and the system may follow either $F_+$ or $F_-$ and return to the switching surface, and may slide yet again. Then it may get to the opposite boundary of the sliding region and follow again $F_+$ or $F_-$, and then yet again return to the switching flow. The exact dynamics would need to be studied further. However, this dynamics differs from Scenario 1 by the fact that the full flow contains sliding evolution.  
\paragraph{Scenario $T_+<S_-<T_-<S_+$ (row 3 in Table~\ref{tab:tab1}, see Figure~\ref{fig_top1})}
In this scenario the flow will not visit $\mathcal{H}_+$ (thus, not follow $F_+$ in the long term). The flow may permanently follow $F_\mathrm{s}$ (sliding), or show infinitely many excursions into $\mathcal{H}_-$ (following $F_-$). The case $S_-<T_+<S_+<T_-$ (row 6 in Table~\ref{tab:tab1}) is a mirror image of this scenario when interchanging subscripts $+$ and $-$. There the flow can have infinitely many excursions into $\mathcal{H}_+$, but will not visit $\mathcal{H}_-$.
%The third scenario we describe is depicted in Fig.~\ref{fig_top1} Case 3 and Case 6. Note that now the sets $T_+$, $T_-$, $S_-$ and $S_+$ follow in a sequence  $S_+ > T_- > S_- > T_+$ or $T_- > S_+ > T_+ > S_-$.
%Between sets corresponding to $T_+$ and $T_-$ there exists an attracting sliding set. Moreover, either, between sets $S_+$ and $T_-$, or  $T_+$ and $S_-$, there exist sets where the vector fields $F_\pm$ point in the same direction with respect to the switching surface. This implies that starting an evolution on the switching surface, we may follow the sliding vector field. If the sliding vector field is stable we are confined to it. Otherwise sets $T_-$ or $T_+$ may be reached, and the system may follow either $F_+$ or $F_-$ and return to the switching surface, and may slide yet again. Then it may reach yet again a boundary of the sliding region and follow again $F_+$ or $F_-$, and then yet again return to the switching flow. The exact dynamics would need to be studied further. The full flow contains sliding evolution and may contain switchings but to one vector field, either to $F_+$ or $F_-$ unlike in Scenario 2 where we may have sliding, then we may follow $F_+$, then sliding again, and then following $F_-$ and so on. Hence, the term sliding with possibly asymmetric switching in our table for this Scenario. 
\paragraph{Scenarios $T_+<S_-<S_+<T_-$ and $T_+<S_+<S_-<T_-$ (rows 4 and 5 in Table~\ref{tab:tab1}, see Figure~\ref{fig_top1})}
%The fourth and the final scenario we describe is depicted in Fig.~\ref{fig_top1} Case 4 and Case 5. Note that now the sets $T_+$, $T_-$, $S_-$ and $S_+$ follow in a sequence $T_- > S_+ > S_- > T_+ $ or $T_- > S_- > S_+ > T_+$.
As the attracting subspaces $\mathcal{M}_{\mathrm{slow},\pm}$ both intersect the switching manifold in the attracting sliding set, the flow permanently follows sliding ($F_\mathrm{s}$).
%This implies that starting an evolution on the switching surface, we may follow the sliding vector field. If the sliding vector field is stable we are confined to it. Otherwise sets $S_-$ or $S_+$ may be reached, and the system may follow $F_s$ along a codimension-two subset.
\section{General return map for repelling sliding case}
\label{sec:swbeh}

In the remaining part of the paper we will focus on the case where the full truncated system has repelling sliding (while the reduced truncated sliding has attracting sliding by assumption). So we have the relations
\begin{align}
    \label{full:rep:slidingFT}
    h_y^0y_{\mathrm{sl},-}<h_x^0f_-^0<h_x^0f_+^0<h_y^0y_{\mathrm{sl},+}.
\end{align}
Since the sliding segment is repelling, initial conditions starting in the sliding segment do not have a unique forward trajectory. However, each point in the repelling sliding segment has admissible non-sliding forward trajectories. For this reason we only consider non-sliding trajectories of the full truncated system \eqref{eq:slfxy} in this section.

We will start with the following Lemma, establishing a return map from the half subspace of the switching subspace
\begin{align*}
    \mathcal{H}_0^+:=\{(x,y)\in\mathcal{H}_0: h_y^0y\geq h_x^0f_+^0\}
\end{align*}
back into itself. The projection $P$ of the set $\mathcal{H}_0^+$ is the part of the vertical axis above point $T_+$ in Figure~\ref{fig_top1}.
\begin{lemma}[Continuous switching for the repelling sliding case]\label{thm:rep:switch}
Let $(x_0,y_0)\in\mathcal{H}_0^+$. There exists a trajectory $(\X_-^t(x_0),\Y_-^t(y_0))$ following $F_-$ and a time $t_1>0$, such that $(\X_-^t(x_0),\Y_-^t(y_0))\in\mathcal{H}_-$ for all $t\in(0,t_1)$, $(x_1,y_1)=(\X_-^{t_1}(x_0),\Y_-^{t_1}(y_0))$ is in the switching subspace $\mathcal{H}_0$ and satisfies $h_y^0y_1\leq h_x^0f_-^0$. 

There also exists a trajectory $(\X_+^t(x_1),\Y_+^t(y_1))$ following $F_+$ and a time $t_2>0$, such that $(\X_+^t(x_1),\Y_+^t(y_1))\in\mathcal{H}_+$ for all $t\in(0,t_2)$, and $(x_2,y_2)=(\X_+^{t_2}(x_1),\Y_+^{t_2}(y_1))$ is back in $\mathcal{H}_0^+$. %Assume initial condition $(x(0),y(0))=(x_0,y_0)\in S_0$ where both vector fields $F_+$ and $F_-$ are transversal to switching surface $S_0$ and point across $S_0$ in the same direction, that is the Lie  derivatives of $h$ in vector fields $F_\pm$ at the initial point $(x_0,y_0)$ have the same sign, which we choose to be negative, that is $\Lie{F_\pm}h(x_0,y_0) < 0$, where $h(x,y) = h_{\mathrm{rd},x}^0x - h_y^0\left[g_y^{-1}\right]^0y$ and the plus/minus symbol indicates states generated by vector fields $F_\pm$. Under the initial assumption, the real part of the spectrum of the eigenvalues of $g_y^0$ is bounded from above by some negative number $-c$. Then there exist times $t_1 > 0$ and $t_2 > 0$ such that $h(x(t),y(t)) < 0$ for $0<t<t_1$, $h(x(t),y(t)) > 0$ for $t_1<t<t_2$ and $h(x(t_1),y(t_1)) = h(x(t_2),y(t_2)) = 0$. We may then express 
An explicit expression for $(x_2,y_2)$ as a function of the initial point $(x_0,y_0)$ and times $t_1$ and $t_2$ is
\begin{align*}
\begin{bmatrix}
    x_2\\y_2
\end{bmatrix}&=\begin{bmatrix}
    x_0 + f_-^0t_1 + f_+^0t_2\\
    \e^{(t_1 + t_2)g_y^0}y_0 + \left[\e^{t_2g_y^0}-\e^{(t_1 + t_2)g_y^0}\right] 
y_{\mathrm{sl},-} + \left[I - \e^{t_2g_y^0}\right]y_{\mathrm{sl},+},
\end{bmatrix}
%y_2 & = .
\end{align*}
where $t_1$ and $t_2$ are functions of $y_0$, but independent of $x_0$.
\label{le:1}
\end{lemma}

\noindent 
{\em Proof.}
The truncated flows for $F_\pm$, given in \eqref{eq:slfxy}, are decoupled in $x$ and $y$, affine and have the explicit formula
\begin{align}
X_\pm^t(x)&=x + f_\pm^0t,&
Y_\pm^t(y)&= y_{\mathrm{sl},\pm} + \e^{tg_y^0}(y - y_{\mathrm{sl},\pm}).
\label{eq:flow}
\end{align}
The switching function of \eqref{eq:slfxy} along such a trajectory is the function  \begin{align}\label{proof:rep:switch:harb}
h_\pm(t;x,y)&=h_{\mathrm{rd},x}^0\left[x + f_\pm^0t\right]+h_y^0\left[g_y^0\right]^{-1}\left[y_{\mathrm{sl},\pm} + \e^{tg_y^0}(y - y_{\mathrm{sl},\pm})\right].
\end{align} For $(x,y)\in\mathcal{H}_0$, we have $h_\pm(0;x,y)=0$ by definition of $\mathcal{H}_0$, such that \eqref{proof:rep:switch:harb} simplifies to
\begin{align}\label{proof:rep:switch:hpm}
h_\pm(t;y)&=h_{\mathrm{rd},x}^0 f_\pm^0t+h_y^0\left[g_y^0\right]^{-1}\left[\e^{tg_y^0}-I\right]\left[y-y_{\mathrm{sl},\pm}\right].
\end{align}
So, $h_\pm$ is independent of $x$. Its time derivative at $t=0$ in some $(x,y)$ and the limit for $t\to\infty$, when following $F_\pm$ are 
\begin{align}\label{proof:rep:switch:hprime0}
\partial_th_\pm(0;y)&=h_x^0f_\pm^0-h_y^0y,&\lim_{t\to\infty}h_\pm(t;y)/t&=h_x^0f_\pm^0-h_y^0y_{\mathrm{sl},\pm},
\end{align}
because $g_y^0$ is a matrix with eigenvalues with negative real parts. 

As $h_y^0y_0>h_x^0f_+^0>h_x^0f_-^0$, vector field $F_-$ points into $\mathcal{H}_-$ at $(x_0,y_0)$, while vector field  $F_+$ points into $\mathcal{H}_-$ or is tangent to $\mathcal{H}_0$, such that we can follow $F_-$.

Using \eqref{proof:rep:switch:hprime0} we observe that $h_-(0;y_0)=0$, $\partial_t h_-(0;y_0)=h_x^0f_-^0-h_y^0y_0<0$ by assumption that $(x_0,y_0)\in\mathcal{H}_0^+$, and that $\lim_{t\to\infty} h_-(t;y_0)/t=h_x^0f_-^0-h_y^0y_{\mathrm{sl},-}>0$ by \eqref{full:rep:slidingFT}.
Thus, $t\mapsto h_-(t;y_0)$ must reach zero for some $t_1>0$. Let us call $(x_1,y_1)=(\X_-^{t_1}(x_0),\Y_-^{t_1}(y_0))$ for this reaching time $t_1$. 

At point $(x_1,y_1)\in\mathcal{H}_0$ we have $\partial_th_-(t_1;y_0)=\partial_th_-(0;y_1)=h_x^0f_-^0-h_y^0y_1\geq0$. Since $h_x^0f_+^0-h_x^0f_-^0$ (see \eqref{full:rep:slidingFT}), we have $h_x^0f_+^0-h_y^0y_1>0$  such that at $(x_1,y_1)$ vector field $F_+$ points into $\mathcal{H}_+$, while $F_-$ points into $\mathcal{H}_+$ or is tangent to $\mathcal{H}_0$, such that we can follow $F_+$. 
We have $h_+(0;y_1)=0$, $\partial_th_+(0;y_1)=h_x^0f_+^0-h_y^0y_1>0$ and $\lim_{t\to\infty}h_+(t;y_1)/t=[h_x^0f_+^0-h_y^0 y_{\mathrm{sl},+}]<0$ according to \eqref{full:rep:slidingFT}. Thus, $t\mapsto h_+(t;y_1)$ must reach zero for some $t_2>0$, and we call $(x_2,y_2)=(\X_+^{t_2}(x_1),\Y_+^{t_2}(y_1))$ for this reaching time $t_2$. The resulting $(x_2,y_2)$ satisfies $(x_2,y_2)\in\mathcal{H}_0^+$ as $\partial_th_+(t_2;y_1)=\partial_th_+(0;y_2)=h_x^0f_+^0-h_y^0y_2\leq0$. 
Combining the expressions for $\X_\pm^t$ and $\Y_\pm^t$, $x_2=\X_+^{t_2}\circ \X_-^{t_1}(x_0)$ and $y_2=\Y_+^{t_2}\circ \Y_-^{t_1}(y_0)$, we express point $(x_2,y_2)$ as
\begin{eqnarray*}
x_2 & = &x_0 + f_-^0t_1 + f_+^0t_2,\\
y_2 & = &\e^{(t_1 + t_2)g_y^0}y_0 + (\e^{t_2g_y^0}-\e^{(t_1 + t_2)g_y^0}) 
y_{\mathrm{sl},-} + (I - \e^{t_2g_y^0})y_{\mathrm{sl},+}.
\end{eqnarray*}

By construction, we had that $t_1$ and $y_1$ are function of $y_0$, but not of $x_0$, such that $t_2$ as a function of $y_1$ (but not of $x_1$) is a function of $y_0$, but not of $x_0$.
\hfill $\square$
%\hfill\qedsymbol

\noindent
Lemma~\ref{thm:rep:switch} permits us to define a return map, which has  a skew-product form (the $x$ component is driven by a map for the $y$ component)
\begin{align*}
    R&:\mathcal{H}_0^+\to\mathcal{H}_0^+,&
    R(x,y)&=(R_x(x,y),R_y(y))\mbox{\ as constructed in Lemma~\ref{thm:rep:switch}.}
\end{align*}
We denote the $x$ and $y$ components of $R$ by $R_x$ and $R_y$ respectively. Furthermore, we denote the travel times constructed as part of the construction in Lemma~\ref{thm:rep:switch} by
\begin{align*}
    t_\pm&:\mathcal{H}_{0,y}^+:=\{y\in\R^m:h_y^0y\geq h_x^0 f_+^0\}\to(0,\infty),\\
    t_-&(y):=t_1,\quad
    t_+(y):=t_2\mbox{\quad as in Lemma~\ref{thm:rep:switch},}
\end{align*}
recalling that (since $h_{\mathrm{rd},x}^0\neq0$) we can find for every $y\in\R^m$ a $x\in\R^n$ such that $(x,y)\in\mathcal{H}_0$. If we restrict $y$ to the half space $\{y\in\R^m:h_y^0y\geq0\}$, then $(x,y)\in\mathcal{H}_0^+$. Furthermore, the travel time will be independent of the $x$ we choose, amking the definition of $t_\pm$ well defined. Using $t_\pm$ we can express $R_x$ as
\begin{align}
  \label{eq:Rxdef}
  R_x(x,y)=x+f_-^0t_-(y)+f_+^0t_+(y).
\end{align}
\begin{corollary}[Eventual boundedness of travel times and solutions]\label{thm:bounds}
There exists a closed bounded subset $\mathcal{H}_{\infty,y}^+\subseteq\mathcal{H}_{0,y}^+$ such that for every $y\in\mathcal{H}_{\infty,y}^+$ the orbit $R_y^\ell(y)$ will eventually stay in $\mathcal{H}_{\infty,y}^+$. A bound for the iterate $\ell_{\min}(y)$ with $R^\ell(y)\in\mathcal{H}_{\infty,y}^+$ for all $\ell\geq\ell_{\min}(y)$ depends only on the norm of $y$, such that we can write $\ell_{\min}(|y|)$.

There exist bounds $0<t_\mathrm{low}<t_\mathrm{up}$ such $t_\pm(y)$ satisfy $t_\mathrm{low}\leq t_\pm(y)\leq t_\mathrm{up}$ for all $y\in\mathcal{H}_{\infty,y}^+$.
\end{corollary}
\noindent
\emph{Proof.} In the construction of $t_\pm(y)$ in Lemma~\ref{thm:rep:switch} the lower bound $t_\mathrm{low}$ is implied by Assumption \eqref{full:rep:sliding} for the scenario we consider. There exists a regular matrix $B\in\R^{m\times m}$ such that $\e^{t g_y^0}$ for all $t\geq t_\mathrm{low}>0$ is a contraction in the norm $|\cdot|_B$ induced by the inner product $y^\tran B^\tran B y$. In this norm the flows $\Y_\pm^t$ map into a bounded ball in the norm $|\cdot|_B$.  A possible choice for the set $\mathcal{H}_{\infty,y}^+$ is the closure of this bounded ball. The time it takes for the flow to map into $\mathcal{H}_{\infty,y}^+$ depends only on the norm of the initial value $y$. For $y\in\mathcal{H}_{\infty,y}^+$ and $t\geq t_\mathrm{low}$ the term $h_y^0\left[g_y^0\right]^{-1}\left[\e^{tg_y^0}-I\right]\left[y-y_{\mathrm{sl},\pm}\right]$ in \eqref{proof:rep:switch:hpm} is bounded, such that the limit for large $t$ is uniform for $y\in\mathcal{H}_{\infty,y}^+$, providing a uniform upper bound $t_\mathrm{up}$ for $t_\pm(y)$.\hfill $\square$

For an orbit $R^\ell(x,y)$ the $x$ component satisfies (denoting the $x$ component of $R^\ell(x,y)$ by $R_x^\ell(x,y)$)
\begin{align}\label{eq:rxdrift}
  \begin{aligned}[c]
    R_x^\ell(x,y)&=x+f_+^0T_{\ell,+}(y)+f_-^0T_{\ell,-}(y)\mbox{,\quad where}\\
    T_{\ell,\pm}(y)&=\sum_{i=0}^{\ell-1}t_\pm(R_y^i(y))\mbox{,\quad defining\ } T_\ell(y)=T_{\ell,+}(y)+T_{\ell,-}(y)
  \end{aligned}
\end{align}
as the orbit's travel time. For $y\in\mathcal{H}_{\infty,y}^+$ Corollary~\ref{thm:bounds} implies that $R_x^\ell(x,y)$ drifts with speed proportional to $\ell$, where the times $T_{\ell,\pm}$ are both bounded by $[\ell t_\mathrm{low},\ell t_\mathrm{up}]$, in the affine plane spanned by $f_\pm^0$, rooted at $x$.

\begin{lemma}[Long-term drift in slow variable]\label{thm:map:xdrift}
Let $(x,y)\in\mathcal{H}_0^+$ be arbitrary. For large $\ell$ the orbit $R^\ell(x,y)$ satisfies
\begin{align}
    \label{eq:map:xdrift}
    \frac{R_x^\ell(x,y)-x}{T_\ell(y)}=f_\mathrm{rd,s}^0+O(1/\ell).
\end{align}
The bounding constant in $O(1/\ell)$ depends only on the norm $|y|$.
\end{lemma}
The vector $f_\mathrm{rd,s}^0$, given in \eqref{eq:slRf3}, is the right-hand side of the truncated reduced sliding flow.

\noindent
\emph{Proof.} Define $\alpha_\ell(y)=T_{\ell,-}(y)/T_\ell(y)$. Then
\begin{align}\label{proof:eq:xdrift}
    \frac{R_x^\ell(x,y)-x}{T_\ell(y)}=f_-^0\alpha_\ell(y)+f_+^0(1-\alpha_\ell(y)).
\end{align}
Since $R^\ell(x,y)$ and $(x,y)$ lie in the switching subspace $\mathcal{H}_0$, they satisfy 
\begin{align*}
    0&=h_{\mathrm{rd},x}^0\left[x+f_+^0T_{\ell,+}(y)+f_-^0T_{\ell,-}(y)\right]-h_y^0\left[g_y^0\right]^{-1}R_y^\ell(y),\\
    0&=h_{\mathrm{rd},x}^0 x-h_y^0\left[g_y^0\right]^{-1}y.
\end{align*}
Subtracting the two identities and dividing them by $T_\ell(y)$, we get
\begin{align*}
    0&=h_{\mathrm{rd},x}^0f_+^0(1-\alpha_\ell(y))+h_{\mathrm{rd},x}^0f_-^0\alpha_\ell(y)-\frac{1}{T_\ell(y)}h_y^0\left[g_y^0\right]^{-1}[R_y^\ell(y)-y].
\end{align*}
We can re-arrange this identity for $\alpha_\ell(y)$, obtaining
\begin{align*}
    \alpha_\ell(y)&=\frac{h_{\mathrm{rd},x}^0f_+^0}{h_{\mathrm{rd},x}^0(f_+^0-f_-^0)} -\frac{1}{T_\ell(y)}\frac{h_y^0\left[g_y^0\right]^{-1}[R_y^\ell(y)-y]}{h_{\mathrm{rd},x}^0(f_+^0-f_-^0)}.
\end{align*}
In the second term the denominator $h_{\mathrm{rd},x}^0f_+^0-h_{\mathrm{rd},x}^0f_-^0$ is a negative (in particular, non-zero) constant, the term $R_y^\ell(y)$ is bounded by $y_\mathrm{bd}=\sup\{|y|:y\in\mathcal{H}_{\infty,y}^+\}$. The denominator $T_\ell(y)$ can be bounded by $T_\ell(y)\geq\ell t_\mathrm{low}$, such that
\begin{align}\label{proof:alpha:est}
    |\alpha_\ell(y)-\alpha_\mathrm{rd}^0|&\leq\frac{1}{\ell}\frac{\|h_y^0\left[g_y^0\right]^{-1}\|(y_\mathrm{bd}+|y|)}{t_\mathrm{low}|h_{\mathrm{rd},x}^0(f_+^0-f_-^0)|}.
\end{align}
This means that replacing $\alpha_\ell(y)$ by $\alpha_\mathrm{rd}^0$ in \eqref{proof:eq:xdrift} results in an error of order $1/\ell$ with a factor as shown on the right in \eqref{proof:alpha:est}, multiplied by $|f_+^0|+|f_-^0|$. This shows the claim of the lemma.\hfill $\square$

The arguments in the proof of Lemma~\ref{thm:rep:switch} also imply the following that from every initial point a trajectory reaches the domain $\mathcal{H}_0^+$ of the map $R$.
\begin{corollary}[Reachability of crossing set]\label{thm:reach}
  For every initial condition $(x,y)$ outside of repelling sliding segment 
  \begin{align*}
      \mathcal{H}_{0,\mathrm{su}}=\{(x,y)\in\mathcal{H}_{0}:h_y^0y\in[h_x^0f_-^0,h_x^0f_+^0]\}
  \end{align*} the trajectory of \eqref{eq:slfxy} will reach $\mathcal{H}_{0}^+$ after a finite time $t_\mathrm{ini}(y)$. For $(x,y)$ inside the repelling sliding segment $\mathcal{H}_{0,\mathrm{su}}$ (where the forward trajectory is not unique) any admissible non-sliding trajectory leaves $\mathcal{H}_{0,\mathrm{s}}$ and also reaches $\mathcal{H}_{0}^+$ after finite time $t_\mathrm{ini}(y)$. The time $t_\mathrm{ini}(y)$ has a uniform upper bound $t_\mathrm{bd}$ for each ball $\{y:|y|\leq y_\mathrm{bd}\}$. 
\end{corollary}
Corollary \ref{thm:reach} and Lemma~\ref{thm:map:xdrift} permit us to make the following statement about the quality of the approximation of a trajectory of the truncated full system \eqref{eq:slfxy} to the sliding trajectory of the truncated reduced system \eqref{eq:slRf}.
\begin{theorem}[Asymptotic sliding approximation]\label{thm:convergence}
  Assume that the truncated full system \eqref{eq:slfxy} satisfies the conditions \eqref{full:rep:slidingFT} for repelling sliding. Let $(x(t),y(t))$ be a trajectory of \eqref{eq:slfxy} starting from $(x(0),y(0))=(x_0,y_0)$ and consisting only of non-sliding segments. Then for large times $t$
  \begin{align}
      \label{eq:convergence}
      \frac{x(t)-x(0)}{t}=f_\mathrm{rd,s}^0+O(1/t).
  \end{align}
  The bound in $O(1/t)$ only depends on $|y|$.
\end{theorem}
Rescaling coordinates back into original $x$, $y$ and time scale results in the statement of Theorem~\ref{thm:main:approx}.
\section{Case $m=1$, one-dimensional fast variable with repelling sliding}
\label{sec:1dfast}

\subsection{Periodic switchings}
\label{sec:dper}
We will show that a natural situation of the dynamics of Scenario  $S_-<T_-<T_+<S_+$ when the fast variable is one-dimensional  implies the existence of a stable fixed point of mapping $R$. This implies further that the equality given in Theorem \ref{thm:convergence} is exact, that is:
$$
 \frac{x(t)-x(0)}{t}=f_\mathrm{rd,s}^0.
$$

 W.l.o.g., we assume that switchings depend on $x^1$ state only, where $x = [x^1,\,x^2,\,\dotsc,\,x^n]$, and also on $y$ state. Hence $h_{\mathrm{rd},x}^0 = [1\,\,0,\dotsc,\, 0]$ and $h_y^0 = 1$. Let $g_y^0 = -a$ ($a > 0$), then the switching surface is given by $\mathcal{H}_0 = \{x^1 + y/a = 0\}$.   
We suppose that we start at some point $(x_0^1,*,y_0)$ where the ``$*$'' symbol indicates remaining $n-2$ initial values of vector field $F_-$. After some $t_1 > 0$ we reach the switching surface again at some point $(x_1^1,*,y_1)$. From $(x_1^1,*,y_1)$ we follow vector field $F_+$ for some $t_2>0$ until the switching surface is reached again at point $(x_2^2,*, y_2)$. Note that this scenario is ensured under conditions stated in Lemma \ref{thm:rep:switch}. We get the following set of equations which define 
the evolution of our system
\begin{eqnarray}
\nonumber
 x(t_1) = x_1 = x_0 + f_-^0t_1,\,\, y(t_1) = y_1 = \mbox{e}^{g_y^0t_1}(y_0 - \left[g_y^{-1}\right]^0g_x^0f_-^0) + \left[g_y^{-1}\right]^0g_x^0f_-^0,\\  \label{eq:C1} \\ 
 \nonumber
 x(t_2) = x_2 = x_1 + f_+^0t_2,\,\, y(t_2) = y_2 = \mbox{e}^{g_y^0t_2}(y_1 - \left[g_y^{-1}\right]^0g_x^0f_+^0) + \left[g_y^{-1}\right]^0g_x^0f_+^0,\\ \label{eq:C2} \\ 
 \nonumber
 h_{\mathrm{rd},x}^0x_0 - \frac{1}{g_y^0}h_y^0y_0 = 0,\,\, h_{\mathrm{rd},x}^0x_1 - \frac{1}{g_y^0}h_y^0y_1 = 0,\,\,h_{\mathrm{rd},x}^0x_2 - \frac{1}{g_y^0}h_y^0y_2 = 0. \\ \label{eq:C3}
\end{eqnarray}
Let us now suppose that $h_{\mathrm{rd},x}^0x_2 = h_{\mathrm{rd},x}^0x_0$ and $y_2 = y_0$. 
It is natural to ask whether such a scenario is possible, and if so, how generic it is, and finally what are the consequences of such a motion. In what follows, such a motion will be termed a period-one orbit. Rearranging and combining equations (\ref{eq:C1})-(\ref{eq:C3}) give the conditions for the existence of such an orbit, which can be expressed as:
\begin{eqnarray}
 & & 0 = ax_0^1+y_0 = 0, \label{eq:p1c1}\\
 & & 0 = ax_0^1+ah_{\mathrm{rd},x}^0 f_-^0t_1+\mbox{e}^{-at_1}(y_0-\left[g_y^{-1}\right]^0g_x^0f_-^0) + \left[g_y^{-1}\right]^0g_x^0f_-^0, \label{eq:p1c2}\\
 & & 0  =  t_1 + t_2\frac{h_{\mathrm{rd},x}^0f_+^0}{h_{\mathrm{rd},x}^0f_-^0},\label{eq:p1c3}\\
 & & y_0 - \left[g_y^{-1}\right]^0g_x^0f_+^0 = \mbox{e}^{-a(t_1 + t_2)}(y_0-\left[g_y^{-1}\right]^0g_x^0f_-^0) + \mbox{e}^{-at_2}\left[g_y^{-1}\right]^0g_x^0(f_-^0-f_+^0). \nonumber \\\label{eq:p1c4}
\end{eqnarray}
Conditions (\ref{eq:p1c1}) and (\ref{eq:p1c2}) imply that points $(x_0,y_0)$ and $(x_1,y_1)$ lie on the switching manifold, whereas conditions (\ref{eq:p1c3}) and (\ref{eq:p1c4}) result from periodicity conditions $x_0^1 = x_2^1$ and $y_0 = y_2$. It can be shown that under the conditions of Lemma \ref{le:1}, matrix $\left[g_y^{-1}\right]^0g_x^0$ may be found which satisfies the above equations 
for $t_1$ ($t_2$) and $y_0$ for an open set of the matrix values and hence a period-1 orbit is a generic case. Now, what are the consequence of such  dynamics? 
Consider equation (\ref{eq:p1c3}), which relates the times of evolution of flows generated by vector fields $F_\pm$.
We have that 
\begin{equation}
\frac{t_1}{t_1+t_2} = -\frac{h_{\mathrm{rd},x}^0f_+^0}{h_{\mathrm{rd},x}^0(f_-^0-f_+^0)} = \alpha,
\label{eq:ftr}
\end{equation}
where $\alpha$ is given by equation (\ref{eq:alpha}), and $\alpha$ represents the same ratio of times as the one in equation (\ref{eq:ftr}) above, but $\alpha$ gives it for the reduced sliding flow. Hence, if there is a stable period-1 solution in the full flow then the dynamics of the full flow corresponds to the dynamics of the sliding flow of the reduced system in that the drift along the slow dynamics is equivalent to the motion of the sliding flow in the reduced system, if we consider projecting the full flow onto the sliding flow. Or putting it yet another way, the equality given in Theorem \ref{thm:convergence} is exact. 
\paragraph{The existence of period-1 cycles}
We focus further on equations (\ref{eq:p1c1})-(\ref{eq:p1c4}). In particular, note that w.l.o.g. we may scale out $a$ parameter as it corresponds to time rescaling. In other words we set $a = 1$. Consider now equation (\ref{eq:p1c3}) which defines time $t_1$ as a function of periodic point $y_0$. We have:
\begin{equation}
-y_0 + h_{\mathrm{rd},x}^0 f_-^0 t_1 + \mbox{e}^{-t_1}(y_0 - \left[g_y^{-1}\right]^0g_x^0f_-^0) + \left[g_y^{-1}\right]^0g_x^0f_-^0 = 0.
\label{eq:py0}
\end{equation}
Let us suppose that $t_1$ is small, so we may expand expression (\ref{eq:py0}) in $t_1$ to leading order. We get:
$$
-y_0 + h_{\mathrm{rd},x}^0f_-^0t_1 + (1-t_1)(y_0-\left[g_y^{-1}\right]^0g_x^0f_-^0) + \left[g_y^{-1}\right]^0g_x^0f_-^0 = 0,
$$
which simplifies to:
\begin{equation}
 t_1(h_{\mathrm{rd},x}^0f_-^0 - y_0 + \left[g_y^{-1}\right]^0g_x^0f_-^0) = 0.
\label{eq:py0s} 
\end{equation}
Note that for (\ref{eq:py0s}) to be equal to $0$, we require 
$$
h_{\mathrm{rd},x}^0f_-^0 - y_0 + \left[g_y^{-1}\right]^0g_x^0f_-^0 = h_x^0f_-^0 - y_0 = 0,
$$
which is satisfied at the boundary of the sliding set of the full system. 
Hence, the existence of period-1 cycles is bounded, for small times $t_1$, by the existence of the 
sliding set in the full system. 

We assume that $h_{\mathrm{rd},x}^0f_-^0 - y_0 + \left[g_y^{-1}\right]^0g_x^0f_-^0 < 0$ (see Theorem~\ref{thm:stable:fixedpoint}) and so, typically, for sufficiently small times, $t_1$ and $t_2$ will be of the order given by the ratio:
\begin{equation}
t_1\approx \frac{y_0 - \left[g_y^{-1}\right]^0g_x^0f_-^0}{h_{\mathrm{rd},x}^0f_-^0}, 
\label{eq:t1approx}
\end{equation}
which implies further that $y_0 - \left[g_y^{-1}\right]^0g_x^0f_-^0 > 0$. Another important value of $y_0$ is $y_0 = \left[g_y^{-1}\right]^0g_x^0f_-^0$. We note that on the set $\{y_0 = \left[g_y^{-1}\right]^0g_x^0f_-^0\}$ we have $\dot y = 0$, and hence this set defines the slow manifold of the full system, which, in turn, defines another boundary for the existence of period-1 cycles.  
\subsection{Stable period-1 cycles}
Clearly, if the system exhibits stable period-1 cycles as defined in the previous section then the singular perturbation of the sliding flow does not significantly alters system dynamics. Hence, it is of our interest to 
determine the set of parameters for which such a scenario takes palace. To this aim consider 
a map from some $(x_0,*,y_0)\in \mathcal{H}_0^+$ back to $\mathcal{H}_0^+$ under the composition of flows  $(\X_+^{t_2}(x_1),\Y_+^{t_2}(y_1))\circ (\X_-^{t_1}(x_0),\Y_-^{t_1}(y_0))$. We drop the superscript ``1'' in the $x_0^1$ component for the clarity of the notation. 
We know that such a map exists under the conditions of Lemma \ref{le:1}. 
We have that 
\begin{equation}
y_2 = \mbox{e}^{-(t_1 + t_2)}(y_0 - \left[g_y^{-1}\right]^0g_x^0f_-^0) + \mbox{e}^{-at_2}\left[g_y^{-1}\right]^0g_x^0(f_-^0-f_+^0)+\left[g_y^{-1}\right]^0g_x^0f_+^0.
\label{eq:mpy}
\end{equation}
Obviously $t_1$ and $t_2$ are functions of $y_0$. We may linearize above map (\ref{eq:mpy}) about some $y_0^*\in \mathcal{H}_0^+$. In the computed derivative below we do not use the star symbol, but $y_0 = y_0^*$, $y_1 = y_1^*(y_0^*)$, $t_1 = t_1(y_0^*)$ and  $t_2 = t_2(y_1^*(y_0^*))$. 
We get:
\begin{eqnarray}
\frac{dy_2}{dy_0} & = &
\mbox{e}^{-(t_1 + t_2)}\left(1 - (y_0 - \left[g_y^{-1}\right]^0g_x^0f_-^0)\frac{1-\mbox{e}^{-t_1}}{h_{\mathrm{rd},x}^0 f_-^0 -\mbox{e}^{-t_1}(y_0-\left[g_y^{-1}\right]^0g_x^0f_-^0)}\right)\nonumber \\
& \times& \left(1 - (y_1 - \left[g_y^{-1}\right]^0g_x^0f_+^0)\frac{1-\mbox{e}^{-t_2}}{h_{\mathrm{rd},x}^0 f_+^0 - \mbox{e}^{-t_2}(y_1-\left[g_y^{-1}\right]^0g_x^0f_+^0)}\right).
\label{eq:mapyd}
\end{eqnarray}
For large enough times $t_i = \mathcal{O}(1)$ ($i = 1,2$), we have 
\begin{equation}
\frac{dy_2}{dy_0}  \approx 
\mbox{e}^{-(t_1 + t_2)}\left(\frac{h_{\mathrm{rd},x}^0f_-^0 - y_0 + \left[g_y^{-1}\right]^0g_x^0f_-^0}{h_{\mathrm{rd},x}^0 f_-^0}\right) \left(\frac{h_{\mathrm{rd},x}^0f_+^0 - y_1 + \left[g_y^{-1}\right]^0g_x^0f_+^0}{h_{\mathrm{rd},x}^0 f_+^0}\right)
\label{eq:mapyda}
\end{equation}
which is small and positive, and tends to $0$ for $t_i \rightarrow \infty$  ($i = 1,\, 2$), by the conditions of Lemma \ref{le:1}. %It can be shown further that for small enough  $t_i$ ($i = 1,\, 2$) 
%$$
%\frac{dy_2}{dy_0} \rightarrow 1.
%$$
%Since the map $y_0 \mapsto y_2$ is strictly monotonic and such that $ 0 < y_2^\prime < 1$ a singular perturbation in our $3D$ systems implies the existence of a globally stable period-1 cycle and hence the singular perturbation in this case may not produce complex dynamics.
\begin{figure}
\centering{\includegraphics[scale=0.4]{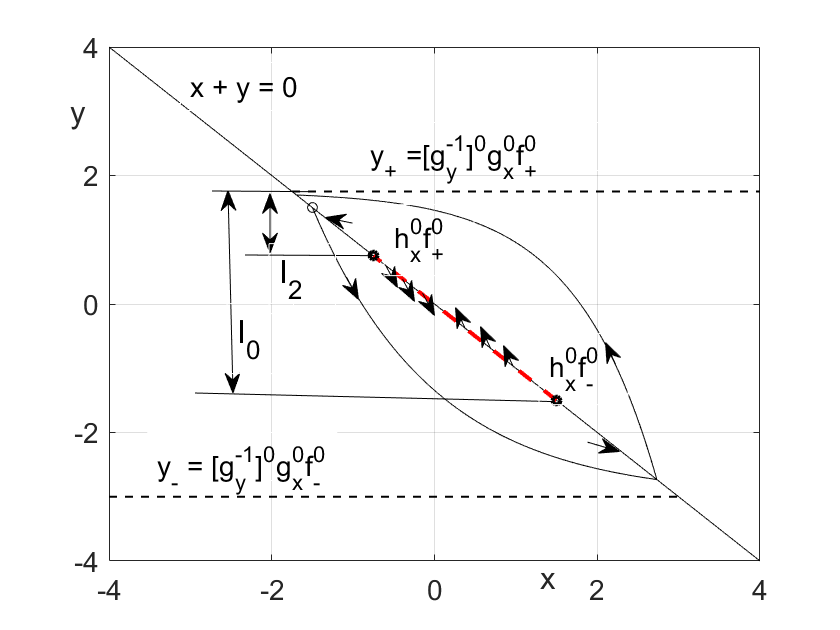}}
\caption{Schematic depiction of phase space projected onto $(x,y)$-plane. The $\dot y = 0$ nullclines, the switching manifold and repelling sliding subset are depicted. Compare with Case 1 in Fig.~\ref{fig_top1}.}
\label{fig3}
\end{figure}

We may also consider what happens for small times $t_1$, and $t_2$. Expanding the exponential terms in the leading order then gives:
\begin{eqnarray}
 \frac{dy_2}{dy_0} & \approx &
(1-(t_1 + t_2))\left(\frac{h_{\mathrm{rd},x}^0f_-^0 - y_0 + \left[g_y^{-1}\right]^0g_x^0f_-^0}{h_{\mathrm{rd},x}^0 f_-^0 - y_0 + \left[g_y^{-1}\right]^0g_x^0f_-^0 + t_1(y_0 - \left[g_y^{-1}\right]^0g_x^0f_-^0)}\right)\nonumber\\
&\times & \left(\frac{h_{\mathrm{rd},x}^0f_+^0 - y_1 + \left[g_y^{-1}\right]^0g_x^0f_+^0}{h_{\mathrm{rd},x}^0 f_+^0 -y_1 + \left[g_y^{-1}\right]^0g_x^0f_+^0 + t_2(y_1- \left[g_y^{-1}\right]^0g_x^0f_+^0)}\right).
\label{eq:mapydasmall}  
\end{eqnarray}
Equivalently, we may express equation (\ref{eq:mapydasmall}) as 
\begin{eqnarray*}
 \frac{dy_2}{dy_0} & \approx &
(1-(t_1 + t_2))\left(\frac{h_x^0f_-^0 - y_0}{h_x^0 f_-^0 - y_0 + t_1(y_0 - \left[g_y^{-1}\right]^0g_x^0f_-^0)}\right)\nonumber\\
&\times & \left(\frac{h_x^0f_+^0 - y_1}{h_x^0 f_+^0 -y_1 + t_2(y_1- \left[g_y^{-1}\right]^0g_x^0f_+^0)}\right).
\end{eqnarray*}

We have that $y_0 - \left[g_y^{-1}\right]^0g_x^0f_-^0 > 0$, $y_1- \left[g_y^{-1}\right]^0g_x^0f_+^0 < 0$, $h_{\mathrm{rd},x}^0 f_-^0 - y_0 + \left[g_y^{-1}\right]^0g_x^0f_-^0 < 0$ and $h_{\mathrm{rd},x}^0f_+^0 - y_1 + \left[g_y^{-1}\right]^0g_x^0f_+^0 > 0$. 
The necessary conditions for times $t_1$ and $t_2$ to be small is that the distances $||(x_1,y_1) - (x_0,y_0)||$, and $||(x_2,y_2) - (x_1,y_1)||$ are  small, which translates on the size of the sliding subset as measured by its euclidian length on the switching surface on the $(x,y)$-plane (since only the first component of the $x$ state is of our interest, we ignore the indices). Moreover due to the uniqueness of the solutions generated by flows $(\X_\pm^{t},\Y_\pm^{t})$, we must have $\frac{dy_2}{dy_0} > 0$, and so it follows that either 
$h_{\mathrm{rd},x}^0 f_-^0 - y_0 + \left[g_y^{-1}\right]^0g_x^0f_-^0 + t_1(y_0 - \left[g_y^{-1}\right]^0g_x^0f_-^0) < 0 \quad \land \quad h_{\mathrm{rd},x}^0 f_+^0 -y_1 + \left[g_y^{-1}\right]^0g_x^0f_+^0 + t_2(y_1- \left[g_y^{-1}\right]^0g_x^0f_+^0) > 0$ or 
$h_{\mathrm{rd},x}^0 f_-^0 - y_0 + \left[g_y^{-1}\right]^0g_x^0f_-^0 + t_1(y_0 - \left[g_y^{-1}\right]^0g_x^0f_-^0) > 0 \quad \land \quad h_{\mathrm{rd},x}^0 f_+^0 -y_1 + \left[g_y^{-1}\right]^0g_x^0f_+^0 + t_2(y_1- \left[g_y^{-1}\right]^0g_x^0f_+^0) < 0$. 
The right-hand side of (\ref{eq:mapydasmall}) is continuous and bounded from the right for $t_i \rightarrow 0$ ($i = 1,2$), and in the limit it is equal to $1$. We will now summarise above findings and make them precise in the theorem below.  
\begin{theorem}\label{thm:stable:fixedpoint}
Let $x\in\mathbb{R}^n$ and $y\in\mathbb{R}$. W.l.o.g. we set $h_{\mathrm{rd},x}^0 = [1,0,\cdots,0]$ which implies that $h_{\mathrm{rd},x}^0x = x^1$. W.l.o.g. we set $a = -1$. Assume the initial condition $(x_0^1,*,y_0)\in \mathcal{H}_0$ where both vector fields $F_+$ and $F_-$ are transversal to $\mathcal{H}_0$ and point across $\mathcal{H}_0$ in the same direction. The ``$*$'' symbol indicates the remaining $n-1$ states of the $x$ state. In this setting, the Lie  derivatives of $h =h_{\mathrm{rd},x}^0+h_y^0y$ in vector fields $F_\pm$ at the initial point $(x_0^1,*,y_0)$ have the same sign, which we choose to be negative, that is $\Lie{F_\pm}h(x_0^1,*,y_0) = -y_0 + h_{\mathrm{rd},x}^0f_\pm^0 + \left[g_y^{-1}\right]^0g_x^0f_\pm^0 < 0$, and the plus/minus symbol indicates states generated by vector fields $F_\pm$. By Lemma \ref{le:1} there exists a return map $y_0\mapsto y_2$, where $y_2$ is given by (\ref{eq:mpy}). 

For $y_0\in (h_x^0f_-^0,\,\,\left[g_y^{-1}\right]^0g_x^0f_+^0]$ the map  (\ref{eq:mpy}) is  monotonically increasing and differentiable from the right in $y_0\downarrow h_x^0f_-^0$ and from the left when $y_0\uparrow \left[g_y^{-1}\right]^0g_x^0f_+^0$, and hence an injection, and it maps interval $I_0 = (h_x^0f_-^0,\,\,\left[g_y^{-1}\right]^0g_x^0f_+^0]$ within $I_2 = (h_x^0f_+^0,\,\,\left[g_y^{-1}\right]^0g_x^0f_+^0)$. For sufficiently large $t_i(y_0)$ ($i = \pm 2$), the derivative $dy_2/dy_0\downarrow 0$, and for sufficiently small $t_i(y_0)$, $dy_2/dy_0\downarrow 1$. It then follows that there exists at least one stable fixed point $y^*\in I_2$ of map  (\ref{eq:mpy}). 
\end{theorem}
{\em Proof}

\begin{figure}
\centering{\includegraphics[scale=0.4]{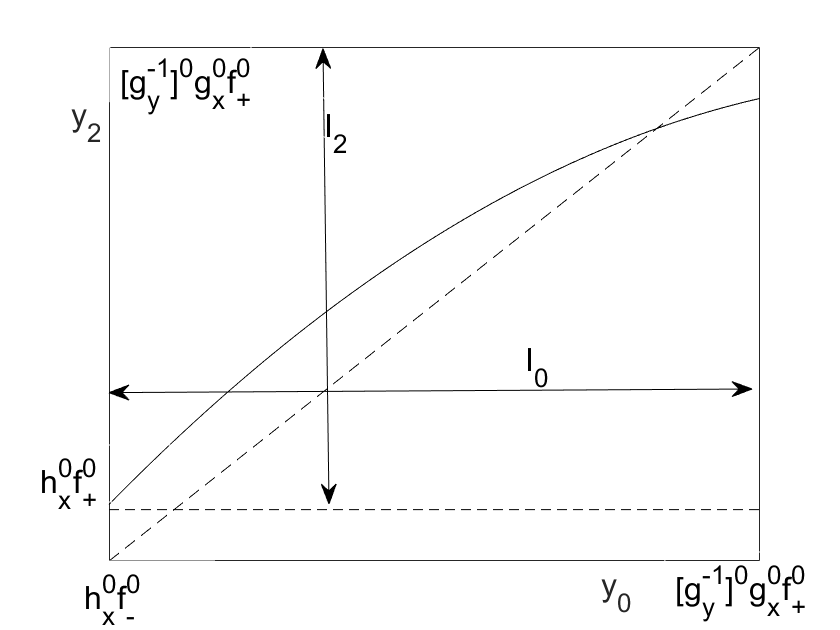}}
\caption{An example of a one-dimensional from $I_0 \ni y_0  \mapsto y_2\in I_2$.}
\label{fig4}
\end{figure}

Starting from any initial $y_0\in I_0$, under the conditions of Lemma \ref{le:1}, the Implicit Function Theorem guarantees that there exist positive times $t_1$ and $t_2$ such that $y_0$ is mapped onto $y_2\in I_2$. Moreover the map (\ref{eq:mpy}) is continuously differentiable for any $y_0\in I_0$ in $y_0$, $t_1(y_0)$ and $t_2(y_0)$ due to the continuity of (\ref{eq:mpy}) and (\ref{eq:mapyd}) with respect to $y_0$, $t_1(y_0)$ and $t_2(y_0)$. Again, to see this we need to invoke the Implicit Function Theorem to determine the existence of continuous derivatives of $t_1(y_0)$ and $t_2(y_0)$ with respect to $y_0$ for any $y_0 \in I_0$. The monotonicity of the map is ensured by the uniqueness of the flow solutions of vector fields $F_\pm$ which make up the two trajectory segments which are concatenated together and map $y_0$ into $y_2$.
By the uniqueness of this concatenation any point starting at some $y_0 \in I_0$ must be mapped to $y_2$, which lies to the left of $y_0\in I_2$, and hence the monotonicity of the map is ensured. Moreover, the former ensures that $dy_2/dy_0 > 0$ for any $y \in I_0$ Finally, we need to show the existence of at least one stable point of map (\ref{eq:mpy}). We note that for any $\vep > 0$ a point $h_x^0f_-^0 + \vep$ must be mapped above point $h_x^0f_+^0$ and hence the image of that point lies above the line $y_2 = y_0$ (see Fig.~\ref{fig4}). Since for sufficiently large $t_i$ ($i = 1,\,2$), which corresponds to points close to $\left[g_y^{-1}\right]^0g_x^0f_+^0$ within $I_0$, the slope $dy_2/dy_0\rightarrow 0$ by equation (\ref{eq:mapyda}), the graph of map (\ref{eq:mpy}) must cross the line $y_2=y_0$ from $y_2 > y_0$ towards $y_2 < y_0$ at some point $y^*\in I_0$. This, in turn, implies that there exists a fixed point $y^*$ in map (\ref{eq:mpy}) such that $ 0 < dy_2/dy_0(y = y^*) < 1$, and hence there exists at least one stable fixed point in (\ref{eq:mpy}). 
It now remains to show that for small $t_i$ ($i = 1,2$), $dy_2/dy_0\rightarrow 1$. Let $y_0 = h_x^0f_-^0+\vep_-$ and 
 $y_1 = h_x^0f_+^0-\vep_+$ where $\vep_\pm$ are sufficiently small and positive. We expand $h(x_1(t_1),y_1(t_1))= 0$ and $h(x_2(t_2),y_2(t_2))= 0$ in $t_1$ and $t_2$, and then solve  
 $$
-y_0 + h_{\mathrm{rd},x}^0f_-^0t_1 + (1-t_1 + \frac{1}{2}t_1^2)(y_0-\left[g_y^{-1}\right]^0g_x^0f_-^0) + \left[g_y^{-1}\right]^0g_x^0f_-^0 = 0,
$$
and 
 $$
-y_1 + h_{\mathrm{rd},x}^0f_+^0t_2 + (1-t_2 + \frac{1}{2}t_2^2)(y_1-\left[g_y^{-1}\right]^0g_x^0f_+^0) + \left[g_y^{-1}\right]^0g_x^0f_+^0 = 0,
$$
for $t_1$ and $t_2$ correspondingly in power series of $\vep_\pm$. To leading order in $\vep_\pm$, 
we find 
\begin{eqnarray}
t_1 & = & \frac{2\vep_-}{h_{\mathrm{rd},x}^0f_-^0}\label{eq:fort1}\\
t_2 & = & -\frac{2\vep_+}{h_{\mathrm{rd},x}^0f_+^0}\label{eq:fort2}.
\end{eqnarray}
Substituting (\ref{eq:fort1}) and (\ref{eq:fort2}) into (\ref{eq:mapydasmall}) gives:
\begin{equation}
\frac{dy_2}{dy_0} = \left(1 - \left(\frac{2\vep_-}{h_{\mathrm{rd},x}^0f_-^0} - \frac{2\vep_+}{h_{\mathrm{rd},x}^0f_+^0}\right)\right)\left(-\frac{1}{1+(2\vep_-)/h_{\mathrm{rd},x}^0f_-^0}\right)\left(-\frac{1}{1-(2\vep_+)/h_{\mathrm{rd},x}^0f_+^0}\right).
\end{equation}
Taking the limit $\vep_\pm\rightarrow 0$ gives:
\begin{equation}
\frac{dy_2}{dy_0} \uparrow 1\quad \mbox{for}\quad \vep_\pm\downarrow 0.
\end{equation}
\hfill $\square$

We illustrate the above theorem in Fig.~\ref{fig3} and Fig.~\ref{fig4}. In particular, in Fig.~\ref{fig3} we depict the phase space projected onto $(x,y)$ components (where $h_{\mathrm{rd},x}x = x^1$ and we drop the upper index 1). In this setting, $x+y = 0$ defines switching manifold $\mathcal{H}_0$, $y_\pm = \left[g_y^{-1}\right]^0g_x^0f_\pm$ are the slow manifolds of vector fields $F_\pm$ respectively. The red dashed line on $\mathcal{H}_0$ denotes the sliding subset within $\mathcal{H}_0$, which is repelling. The boundaries of this set are also represented by $h_x^0f_\pm^0$. Compare with Scenario $S_+ > T_+ > T_- > S_-$ in Fig.~\ref{fig_top1}. 
An example of a trajectory starting at some $y_0$ within $I_0$ is shown to lead to some point $y_2\in I_2$ after a switch at $S_0$. A typical map $I_0 \ni y_0  \mapsto y_2\in I_2$ is then schematically depicted in Fig.~\ref{fig4}. A stable fixed point corresponding to a period-1 cycle is clearly visible as the point of intersection between the curve representing the map with the line $y=x$. 
\section{Numerical examples}
\label{sec:examples}

%\subsection{Switchings and drift along the sliding flow}
We will now present numerical examples which illustrate our findings. We expect that, typically, a sliding flow will be perturbed to a periodic behaviour with a drift. 
\paragraph{Example 1}
We set $h_{\mathrm{rd},x}^0 = [1\,\,\, 0]$, $f_-^0 = [1.5\,\,\, 1]^T$, $f_+^0 = [-1\,\,\, -0.5]^T$, $\left[g_y^{-1}\right]^0g_x^0 = [-1\,\,\, -1.5]$, $-g_y^0 = a = 2$ and $h_y^0 = 1$. The initial conditions are $(x_1(0),x_2(0),y(0)) = (-0.25,0,0.5)$ (full system), $(x_1(0),x_2(0)) = (0,0)$ (reduced system).
 \begin{figure}
\centering{\includegraphics[scale=0.25]{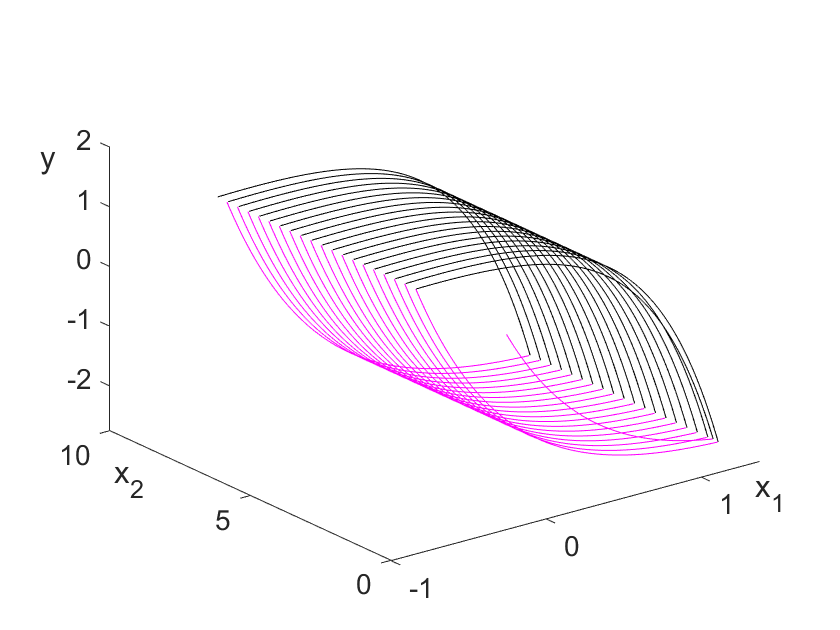}}
\caption{Representative trajectory in phase space of the full system (\ref{eq:slfxy}) for numerical values set to: $h_{\mathrm{rd},x}^0 = [1\,\, 0]$, $f_+^0 = [-1\,\, -1.5]^T$, $f_-^0 = [1.5\,\, 1]^T$, $\left[g_y^{-1}\right]^0g_x^0 = [-1\,\, -1.5]$, $-g_y^0 = a =  2$, $h_y^0 = 1$, and the initial conditions $(x_1(0),x_2(0),y(0)) = (-0.25,0,0.5)$ (full system), $(x_1(0),x_2(0)) = (0,0)$ (reduced system - not depicted).}
\label{fig1}
\end{figure}
A representative trajectory in phase space of the 3-dimensional slow-fast system is depicted in Fig.~\ref{fig1}. We find that $t_1$ and $t_2$ converge to values $t_1 = 1.4857$ and $t_2 = 2.2285$ and hence define a stable period-one orbit for which we find:
$$
\alpha = -\frac{h_{\mathrm{rd},x}^0f_+^0}{h_{\mathrm{rd},x}(f_-^0-f_+^0)} = \frac{t_1}{t_1+t_2} = 0.4.
$$
\paragraph{Example 2}
We now set $h_{\mathrm{rd},x}^0 = [1\,\, 0.5]$, $f_+^0 = [-1\,\, -0.5]^T$, $f_-^0 = [1.5\,\, 1]^T$, $\left[g_y^{-1}\right]^0g_x^0 = [-1\,\, -0.5]$, $-g_y^0 = 2$, $h_y^0 = 1$ and the initial condition $(x_1(0),x_2(0),y(0)) = (-0.25,0,0.5)$ (full system), $(x_1(0),x_2(0)) = (0,0)$ (reduced system). Representative trajectory in phase space in Fig.~\ref{fig2}. Similarly as in Example 1 the initial conditions for the reduced as well as the full system differ in that in each case we want to commence the evolution on the switching manifold. So, one may consider the initial condition of the full system as a perturbation of the initial conditions of the reduced flow. However, in each case the trajectory starts on the switching manifold (in the reduced case we are confined to it).      
\begin{figure}
\centering{\includegraphics[scale=0.25]{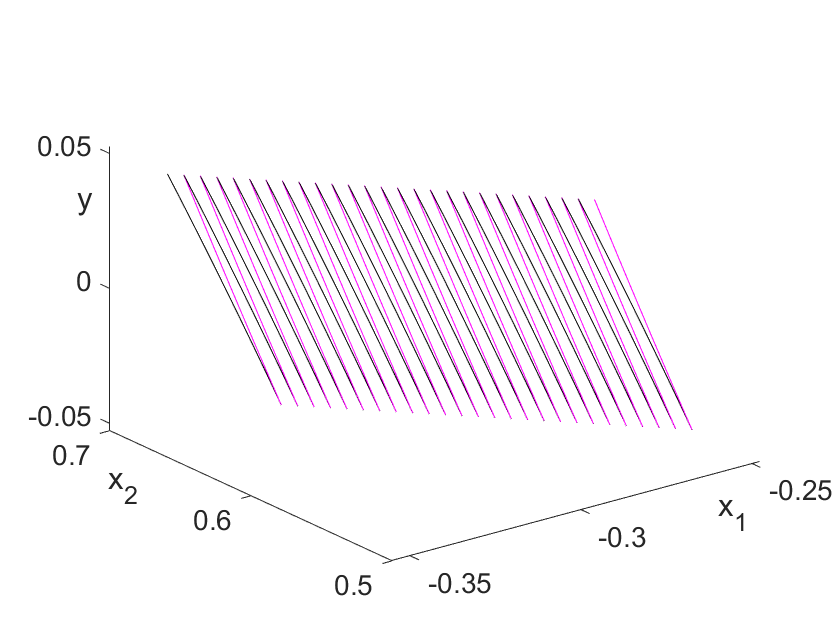}}
\caption{Trajectory in phase space of the full systems (\ref{eq:slfxy}) for numerical values set to: $h_{\mathrm{rd},x}^0 = [1\,\, 0.5]$, $f_+^0 = [-1\,\, -0.5]^T$, $f_-^0 = [1.5\,\, 1]^T$, $\left[g_y^{-1}\right]^0g_x^0 = [-1\,\, -0.5]$, $-g_y^0 = 2$, $h_y^0 = 1$, and the initial conditions $(x_1(0),x_2(0),y(0)) = (-0.25,0,0.5)$ (full system), $(x_1(0),x_2(0)) = (0,0)$ (reduced system - not depicted).}
\label{fig2}
\end{figure}
 We find that $t_1$ and $t_2$ converge to values $t_1 = 0.0263$ and $t_2 = 0.0421$ and hence define a stable period-one orbit for which we find:
$$
\alpha = -\frac{h_{\mathrm{rd},x}^0f_+^0}{h_{\mathrm{rd},x}^0(f_-^0-f_+^0)} = \frac{t_1}{t_1+t_2} = 0.3846.
$$
%\begin{figure}
%\centering{\includegraphics[scale=0.4]{fig2_case1.png}}
%\caption{Representative solutions of reduced (\ref{eq:rS2}) and full (\ref{eq:slfxy}) systems in phase space for numerical values set to: $h_x = [1\,\, 0.5]$, $f_+ = [-1\,\, -0.5]^T$, $f_- = [1.5\,\, 1]^T$, $B = [-1\,\, -0.5]$, $A = 2$, $h_y = 1$, and the initial conditions $(x_1(0),x_2(0),y(0)) = (-0.25,0,0.5)$ (full system), $(x_1(0),x_2(0)) = (0,0)$ (reduced system).}
%\label{fig2}
%\end{figure}
It may be further verified that the sliding subset $\mathcal{H}_{0,\mathrm{su}}$ in Example 2 is smaller than the one in Example 1 (using as measure its euclidian length on $\mathcal{H}_0$ confined to $(x,y)$) as it is implicitly stated in Theorem \ref{thm:stable:fixedpoint}.
\section{Conclusions}
\label{sec:conclusions}
In the current work we investigate the effects of singular perturbation on the sliding dynamics present in a Filippov system. In particular, we consider the case when the sliding flow corresponds to the dynamics of the reduced Filippov system. We then consider what is the effect of stable singular perturbations on that flow. Generically, for $n$-dimensional slow dynamics and $m$-dimensional fast dynamics, due to singular perturbations we distinguish five qualitatively distinct scenarios as we classify them in Sec.~\ref{sec:results} (see Tab.~1 and Fig.~1) by means of a sequence of trajectory segments generated by the perturbed flow. This is our first main result. Namely, due to singular perturbations of the sliding reduced flow we may observe: a) switchings between vecor fields $ F_\pm$ along the sliding flow of the reduced system (Scenario $S_-<T_-<T_+<S_+$); b) sliding or sliding segments possibly interspersed with switchings to $F_\pm$ (Scenario $S_-<T_+<T_-<S_+$); c) sliding or sliding segments possibly interspersed with switchings to $F_-$ (Scenario $T_+<S_-<T_-<S_+$); d) sliding or sliding segments possibly interspersed with switchings to $F_+$ (Scenario $S_-<T_+<S_+<T_-$); e) sliding (Scenarios: $T_+<S_-<S_+<T_-$, $T_+<S_+<S_-<T_-$). In Scenarios b)-d) the exact dynamics will depend on the stability properties of the sliding flow in the full system. We may further classify scenarios a) to e) into two cases. Namely, case one which corresponds to scenario a) which is characterised by the presence of the repelling sliding set in the perturbed systems, and case two which corresponds to scenarios b) to e) which is characterised by the presence of an attracting sliding set in the perturbed system.     

 In the current paper case one is then investigated in detail in Sec. ~\ref{sec:swbeh} and \ref{sec:1dfast}. 
 The main results of Sec. ~\ref{sec:swbeh}, encapsulated in the summary as Theorem~\ref{thm:main:approx},
  states that stable singular perturbations of the sliding flow with the fast dynamics having dimension $m > 1$ results in perturbation of $\mathcal{O}(\vep)$ that depends on $t$ and is only continuous, so it can be rapidly fluctuating. This, in turn, may imply the existence of micro chaotic dynamics as reported in different contexts in Filippov type systems \cite{haller1996micro,csernak2010digital,GlKo2010}. The special case of $m= 1$ is then treated in Sec.~\ref{sec:1dfast}. In such case, stable singular perturbation implies regular perturbation of $\mathcal{O}(\vep)$ of the full flow, and hence no complex dynamics, even on a micro scale, is possible. We illustrate our findings by two numerical examples in Sec. \ref{sec:examples}. The future work is focused on investigating case two and applying our analytical findings to models of relevance to applications. 

To summarise: case one may be considered as ``safe'' in the sense that the singular perturbation may not significantly alter system dynamics. On the other hand, case two may turn out to be ``unsafe'' in the sense that, for example, a significant change in the direction of the perturbed flow cannot be excluded. Of course, it remains to be investigated if there is such a possibility. 
%%%%%%%%%%%%%%%%%%%%%%%%%%%%%%%%%%%%%%%%%%%%%%%%%%%%%%%%%%%%%%%%%%%
\bibliographystyle{spmpsci}
\bibliography{refs}

\begin{thebibliography}{10}
\providecommand{\url}[1]{{#1}}
\providecommand{\urlprefix}{URL }
\expandafter\ifx\csname urlstyle\endcsname\relax
  \providecommand{\doi}[1]{DOI~\discretionary{}{}{}#1}\else
  \providecommand{\doi}{DOI~\discretionary{}{}{}\begingroup
  \urlstyle{rm}\Url}\fi

\bibitem{AdDaNo:01}
Adolfsson, J., Dankowicz, H., Nordmark, A.: {3D} passive walkers: {F}inding
  periodic gaits in the presence of discontinuities.
\newblock Nonlinear Dynamics pp. 205--229 (2001)

\bibitem{diBeBuChaKo:08}
di~Bernardo, M., Budd, C., Champneys, A., Kowalczyk, P.: Piecewise-smooth
  Dynamical Systems: Theory and Applications.
\newblock Springer-Verlag (2008)

\bibitem{COLOMBO20131}
Colombo, A., Jeffrey, M.R.: The two-fold singularity of nonsmooth flows:
  Leading order dynamics in n-dimensions.
\newblock Physica D: Nonlinear Phenomena \textbf{263}, 1--10 (2013)

\bibitem{csernak2010digital}
Csern{\'a}k, G., St{\'e}p{\'a}n, G.: Digital control as source of chaotic
  behavior.
\newblock International Journal of Bifurcation and Chaos \textbf{20}(05),
  1365--1378 (2010)

\bibitem{fe:79}
Fenichel, N.: Geometric singular perturbation theory for ordinary differential
  equations.
\newblock Journal of Differential Equations \textbf{31}(1), 53--98 (1979)

\bibitem{Fi:88}
Filippov, A.: Differential Equations with Discontinuous Righthand Sides.
\newblock Kluwer Academic Publishers, Dortrecht (1988)

\bibitem{Fr:02a}
Fridman, L.: Singularly perturbed analysis of chattering in relay control
  systems.
\newblock IEEE Trans. on Automatic Control \textbf{47}(12), 2079--2084 (2002)

\bibitem{Fr:02b}
Fridman, L.: Slow periodic motions with internal sliding modes in variable
  structure systems.
\newblock Int. J. of Control \textbf{75}(7), 524--537 (2002)

\bibitem{GlKo2010}
Glendinning, P., Kowalczyk, P.: Micro-chaotic dynamics due to digital sampling
  in hybrid systems of filippov type.
\newblock Physica D \textbf{239}(1-2), 44--57 (2010)

\bibitem{haller1996micro}
Haller, G., St{\'e}p{\'a}n, G.: Micro-chaos in digital control.
\newblock Journal of Nonlinear Science \textbf{6}(5), 415--448 (1996)

\bibitem{KaKr:19}
Kaklamanos, P., Kristiansen, K.: Regularization and geometry of piecewise
  smooth systems with intersecting discontinuity sets.
\newblock SIAM Journal on Applied Dynamical Systems \textbf{18}(3), 1225--1264
  (2019)

\bibitem{KrHo:15a}
Kristiansen, K., Hogan, J.: On the use of blowup to study regularizations of
  singularities of piecewise smooth dynamical systems in $\mathbb{R}^3$.
\newblock SIAM Journal on Applied Dynamical Systems \textbf{14}(1), 382--422
  (2015)

\bibitem{KrHo:15b}
Kristiansen, K., Hogan, J.: Regularizations of two-fold bifurcations in planar
  piecewise smooth systems using blowup.
\newblock SIAM Journal on Applied Dynamical Systems \textbf{14}(4), 1731--1786
  (2015)

\bibitem{KrSz:01}
Krupa, M., Szmolyan, P.: Relaxation oscillation and canard explosion.
\newblock Journal of Differential Equations \textbf{174}, 312--368 (2001)

\bibitem{LeNi:04}
Leine, R., Nijmeijer, H.: Dynamics and Bifurcations of Non-Smooth Mechanical
  Systems.
\newblock Lecture Notes in Applied and Computational Mechanics.
  Springer--Verlag, Berlin Heidelberg (2004)

\bibitem{MaLa:12}
Makarenkov, O., Lamb, J.: Dynamics and bifurcations of nonsmooth systems: A
  survey.
\newblock Physica D \textbf{241}(22), 1826--1844 (2012)

\bibitem{SiKo:10}
Sieber, J., Kowalczyk, P.: Small-scale instabilities in dynamical systems with
  sliding.
\newblock Physica D \textbf{239}, 44--57 (2010)

\bibitem{Si:25}
Simpson, D.: Nonsmooth folds as tipping points.
\newblock Chaos \textbf{35}(2) (2025)

\bibitem{SiGl:24}
Simpson, D., Glendinning, P.: Inclusion of higher-order terms in the
  border-collision normal form: Persistence of chaos and applications to power
  converters.
\newblock Physica D \textbf{462} (2024)

\bibitem{Ut:92}
Utkin, V.: Sliding Modes in Control Optimization.
\newblock Springer--Verlag, New York (1992)

\end{thebibliography}

\end{document}